\documentclass[11pt]{article}%
\usepackage{amsfonts}
\usepackage{amsmath}
\usepackage{amssymb}
\usepackage[onehalfspacing]{setspace}
\usepackage{graphicx}%
\newenvironment{proof}[1][Proof]{\noindent\textbf{#1.} }{\ \rule{0.5em}{0.5em}}
\newtheorem{theorem}{Theorem}[section]

\newtheorem{definition}{Definition}[section]

\newtheorem{lemma}{Lemma}[section]
\newtheorem{proposition}{Proposition}[section]
\newtheorem{remark}{Remark}[section]
\begin{document}

\title{Weak convergence in the functional autoregressive model.}
\author{Andr\'{e} Mas\thanks{Institut de Mod\'{e}lisation Math\'{e}matiques de
Montpellier, UMR 5149, Universit\'{e} Montpellier 2, CC051, Place Eug\`{e}ne
Bataillon, 34095 Montpellier Cedex 5, France. (Phone : +33 467143956, Fax :
+33 467143558)\newline mas@math.univ-montp2.fr}}
\date{}
\maketitle

\begin{abstract}
The functional autoregressive model is a Markov model taylored for data of
functional nature. It revealed fruitful when attempting to model samples of
dependent random curves and has been widely studied along the past few years.
This article aims at completing the theoretical study of the model by
adressing the crucial issue of weak convergence for estimates from the model.
The main difficulties stem from an underlying inverse problem as well as from
dependence between the data. Traditional facts about weak convergence in non
parametric models appear : the normalizing sequence is not an $o\left(
\sqrt{n}\right)  $, a bias terms appears. Several original features of the
functional framework are pointed out.

\end{abstract}

\bigskip


\section{Introduction}

\subsection{The model and its history}

The Functional Autoregressive Model of order 1 (FAR1) generalizes to random
elements with values in an infinite dimensional space the classical AR(1)
model belonging to the celebrated class of ARMA process, widely used in time
series analysis. This model was introduced by Bosq \cite{Bos1}, then studied
by several authors. Several chapters in Bosq \cite{Bos3} are dedicated to a
thorough study of this strictly stationary process $\left(  X_{n}\right)
_{n\in\mathbb{Z}}$ defined by%
\begin{equation}
X_{n}-m=\rho\left(  X_{n-1}-m\right)  +\varepsilon_{n,}\quad n\in\mathbb{Z},
\label{modele}%
\end{equation}
where the $X_{k}$'s and the $\varepsilon_{k}$'s are random elements with
values in an infinite dimensional vector space $\mathcal{E}$, $\rho$ is an
unknown linear operator from $\mathcal{E}$ to $\mathcal{E}$ and $m\in
\mathcal{E}$ is the expectation of the process. In all the following we will
assume that for all $n$ $\varepsilon_{n}$ is independent of $X_{n-1}$. The
process $\left(  X_{n}\right)  _{n\in\mathbb{Z}}$ is Markov whenever the
$\varepsilon_{n}$'s are such that $\mathbb{E}\left(  \varepsilon_{n}%
|X_{n-1}\right)  =0$ where $\mathbb{E}$ denotes expectation.

The model was extended in Mourid \cite{Mou} considering autoregressive
processes of higher orders. Besse and Cardot \cite{BC} proved that the model
is adapted to splines techniques. Then Pumo \cite{Pu} studied autoregressive
processes with values in the Banach space of continuous functions on $\left[
0,1\right]  .$ The PhD\ Thesis by Mas \cite{M00} was partly devoted to the
topic. Besse, Cardot and Stephenson \cite{BCS} developped a method based on
kernels. Recently Mas and Menneteau \cite{M3} announced large and moderate
deviations theorems for the process or its covariance sequence whereas
Antoniadis and Sapatinas \cite{Anto} implemented wavelet methods which
considerably improved the prevision mean square error. Even more recently
Menneteau \cite{Men} proved laws of the iterated logarithm for statistics
arising from functional PCA of the process.

The model revealed fruitful in several areas of applied statistics :
electrical engineering (Cavallini et alii, \cite{CMLLC}), climatology (Besse
et alii, \cite{BCS}, Antoniadis and Sapatinas \cite{Anto}), medicine (Marion
and Pumo, \cite{MaPu}).

The main interest of (\ref{modele}) relies in its predictive power. Estimating
the correlation operator only aims at providing an estimate, say $\rho_{n}$
yielding a predictor for the unknown $X_{n+1},$ $\rho_{n}(X_{n})$ based on the
sample $\left(  X_{1},...,X_{n}\right)  $.

However if convergence of $\rho_{n}(X_{n})$ to $\rho(X_{n})$ for instance was
often studied either in probability or almost surely, the issue of weak
convergence has not been truly tackled yet. An attempt was proposed in Mas
\cite{M99} but the conditions under which the result holds are extremely
restricting. The problem of weak convergence is especially intricate due to
the functional framework and to an underlying inverse problem (see next
section). A weak convergence result implies obtaining the sharpest rate for
convergence in probability. Authors studying rates of convergence for the
predictor usually just give bounds... Besides a weak convergence result would
be of much help in getting confidence sets for $\rho(X_{n})=E\left(
X_{n+1}|X_{n},X_{n-1},...\right)  $. Maybe a bootstrap procedure could be
proposed to achieve the same goal but on a one hand I did not find any real
and reliable bootstrap procedure adapted to this pure functional framework in
the literature. On the other hand even if a bootstrap approach may be
satisfactory on a practical viewpoint, it will just provide an approximate
distribution. Here the \textit{exact} asymptotic distribution is given.
Besides the scope of the paper is rather theoretical. The surprising Theorem
\ref{TH1} for instance is -to me at last- really food for thought for people
dealing with functional data. However a promising approach would be to compare
the results of this paper and those obtained by a bootstrap procedure, if any
is available.

One of the other interests of the model is its simplicity. However in the
general framework mentioned above, a first problem arises : in the case of a
general space $\mathcal{E}$, not much is known about the mathematical
description and properties of the linear space, say $\mathcal{L}\left(
\mathcal{E}\right)  $, of bounded linear operators from $\mathcal{E}$ to
$\mathcal{E}$. Estimating $\rho$ requires to build a sequence $\rho_{n}$ of
random linear operators in $\mathcal{L}\left(  \mathcal{E}\right)  ,$ and we
may face serious troubles if the space $\mathcal{L}\left(  \mathcal{E}\right)
$ is too complex.

Usually authors focus on special cases and take for instance $\mathcal{E=C}%
^{m}\left[  0,T\right]  ,$ a Banach space of functions defined on $\left[
0,T\right]  $ and with several continuous derivatives (in Mourid, \cite{Mou})
or $\mathcal{E=}W^{m,p}\left[  0,T\right]  ,$ a space of Sobolev functions on
a real interval (see for instance Adams \cite{Adams}) for definitions and
properties of Sobolev spaces). There are practical reasons for these choices.
Indeed, the curves $X_{n}$ are observed at discretized times and must be first
reconstructed by implementing splines or wavelets for instance. These
techniques provide explicit functions belonging to the spaces mentioned above.

Here appears the second problem : studying weak convergence for random
elements, such as our predictor $\rho_{n}(X_{n})$, in general infinite
dimensional spaces is especially difficult, sometimes tricky. The most general
tool is the Portmanteau Theorem (see Billingsley, \cite{Bil}) but it is rather
a general definition than a criterion to check the convergence of measures.
Even if we consider the Central Limit Theorem which is a very important but
special case of convergence in distribution for measures, there are only a few
spaces for which sufficient conditions are available (even fewer for a
necessary condition). We refer to Ledoux and Talagrand \cite{LT} for a review
on the CLT in Banach spaces. However, if $\mathcal{E}$ is a \textbf{separable
Hilbert space}, the situation becomes more favourable. Take $Z_{i}$ a sequence
of random elements in $\mathcal{E}$. It is a well-known fact that the CLT
holds for i.i.d. $Z_{i}$ if and only if the strong second moment is finite
(i.e. $\mathbb{E}\left\Vert Z_{i}\right\Vert ^{2}<+\infty$). Besides many
authors studied the CLT under different sorts of dependence assumptions
(m-dependence, mixing, martingale differences, etc). We refer to Araujo-Gine
\cite{AG} for a monograph on the CLT. The Hilbertian setting is quite
comfortable for several other well-known mathematical reasons :

\begin{itemize}
\item All Hilbert spaces are isometrically isomorphic to the sequence space
$l^{2}$ , hence have the same underlying geometric structure. They appear as
the most natural generalization of the Euclidean space to the infinite
dimensional setting.

\item The bases are denumerable, the paralellogram identity is valid, the
projection on convex sets is uniquely defined.

\item The operator $\rho$ belongs to the Banach space of linear operators on a
Hilbert space.\ This space is widely used in several areas of mathematics.
Spectral decompositions are available for compact operators.
\end{itemize}

In all the sequel we will set once and for all $\mathcal{E}=\mathcal{H}$ and
$\mathcal{H}$ will usually be a space $W^{m,2}$ where the smoothness index $m$
belongs to $\mathbb{N}$ ($W^{0,2}=L^{2}$).

The next remark is related to $\rho$ and also aims at restricting the field of
our research in order to gain some accuracy in the forthcoming results. In
fact the space $\mathcal{L}\left(  \mathcal{H}\right)  $ is much too large :
this Banach space is not separable. This could turn out to be a serious
problem as far as measurability is concerned (remind that we need to define a
sequence of estimates $\rho_{n}$ for $\rho$ taking values in $\mathcal{L}%
\left(  \mathcal{H}\right)  $). For other reasons mentioned in the next
section, we will suppose that $\rho$ is a \textbf{compact operator}. The space
$\mathcal{K}\left(  \mathcal{H}\right)  $ of compact operators is separable,
its properties are closed to those of (finite size) matrices. Many features of
linear operators on finite dimensional spaces are generalized to
$\mathcal{K}\left(  \mathcal{H}\right)  $ in a kind way.

The space $\mathcal{H}$ is endowed with norm $\left\Vert \cdot\right\Vert ,$
derived from the scalar product $\left\langle \cdot,\cdot\right\rangle .$ In
the case where $\mathcal{H}=W^{m,2}$ we have
\[
\left\langle f,g\right\rangle =\sum_{j=0}^{m}\int f^{\left(  j\right)
}\left(  s\right)  g^{\left(  j\right)  }\left(  s\right)  ds.
\]
Spaces of continuous operators on $\mathcal{H}$ are endowed with the classical
sup-norm defined for all bounded operator $T$ by
\[
\left\Vert T\right\Vert _{\infty}=\sup_{x\in\mathcal{H}_{1}}\left\Vert
Tx\right\Vert
\]
where $\mathcal{H}_{1}$ is the unit ball of $\mathcal{H}$.

The space of Hilbert-Schmidt operators denoted $\mathcal{K}_{2}\left(
\mathcal{H}\right)  $ is endowed with norm $\left\Vert T\right\Vert _{2}%
=\sum_{p}\left\Vert T\left(  e_{p}\right)  \right\Vert ^{2}$ where $e_{p}$ is
any c.o.n.s. in $\mathcal{H}$. The spaces $\mathcal{K}_{2}\left(
\mathcal{H}\right)  $ is a subspace of $\mathcal{K}\left(  \mathcal{H}\right)
$. Note that up to the author's knowledge, the literature on model
(\ref{modele}) or its close alternatives in an Hilbertian framework assumes
that $\rho\in\mathcal{K}_{2}\left(  \mathcal{H}\right)  $. Consequently we
consider in this article a larger class for the unknown parameter.

The tensor product notation is of much use. It enables to define finite rank
operators. For $u,v\in\mathcal{H},$
\[
\left(  u\otimes v\right)  \left(  h\right)  =\left\langle u,h\right\rangle
v.
\]
We may have to deal with another space of operators : the space of trace class
operators $\mathcal{K}_{1}\left(  \mathcal{H}\right)  \subset\mathcal{K}%
\left(  \mathcal{H}\right)  \mathcal{\ }$(the $\left\Vert \cdot\right\Vert
_{1}$ norm on this space will not be fully defined here but I just mention
that $\left\Vert u\otimes v\right\Vert _{1}=\left\Vert u\right\Vert \left\Vert
v\right\Vert $). Finally we will sometimes use the following norm bound :%
\[
\left\Vert \cdot\right\Vert _{\infty}\leq\left\Vert \cdot\right\Vert _{2}%
\leq\left\Vert \cdot\right\Vert _{1}.
\]

\section{Identification and covariance regularization}

In this Hilbert space setting, Bosq \cite{Bos3} proved that whenever it exists
$j_{0}$ such that $\left\Vert \rho^{j_{0}}\right\Vert _{\infty}<1$ and when
$\mathbb{E}\left\Vert \varepsilon_{1}\right\Vert ^{2}$ is finite, $X_{n}$ is a
strictly stationary sequence. For the sake of simplicity and in order to
alleviate calculations within the proofs we will assume that $\left\Vert
\rho\right\Vert _{\infty}<1$. In the sequel we will assume that $\mathbb{E}%
\left(  X_{n}\right)  =m=0$ i.e. we will not adress the problem of estimating
the mean since this issue was extensively treated in the literature. But we
have to face two other serious issues.

\subsection{Identifiability}

As the data are of functional nature, the inference on $\rho$ cannot be based
on likelihood. Lebesgue's measure does not exist on non locally compact spaces
and up to the author's knowledge the classical notion of density has not been
extended to functional random elements. A classical \textbf{moment method}
provides the following normal equation :%
\begin{equation}
\Delta=\rho\Gamma\label{mom}%
\end{equation}
where
\begin{align*}
\Gamma &  =\mathbb{E}\left(  X_{1}\otimes X_{1}\right)  ,\\
\Delta &  =\mathbb{E}\left(  X_{1}\otimes X_{2}\right)
\end{align*}
are the covariance operator (resp. the cross covariance operator of order one)
of the process $\left(  X_{n}\right)  _{n\in\mathbb{Z}}$.

It is a well known fact that whenever $\mathbb{E}\left(  \left\Vert
X_{1}\right\Vert ^{2}\right)  $ is finite $\Gamma$ is a selfadjoint positive,
trace class operator (hence compact). In other words, $\Gamma$ admits the
following Schmidt (i.e. spectral) decomposition :%
\begin{equation}
\Gamma=\sum_{l\in\mathbb{N}}\lambda_{l}\pi_{l},\quad\sum_{l\in\mathbb{N}%
}\lambda_{l}<+\infty\label{schmidt}%
\end{equation}
where $\left(  \lambda_{l}\right)  _{l\geq1}$ is the sequence of the positive
eigenvalues of $\Gamma$ and $\left(  \pi_{l}\right)  _{l\geq1}$ is the
associated sequence of projectors. In the sequel the eigenvectors of $\Gamma$
are denoted $\left(  e_{l}\right)  _{l\geq1}$ hence $\pi_{l}=e_{l}\otimes
e_{l}$ and if $x$ is any vector of $\mathcal{H}$ we set $x_{p}=\left\langle
x,e_{p}\right\rangle .$ For further purpose $\Gamma_{\varepsilon}%
=\mathbb{E}\left(  \varepsilon_{1}\otimes\varepsilon_{1}\right)  $ will stand
for the covariance operator of $\varepsilon_{1}$.

The first step consists in checking that equation (\ref{mom}) correctly
defines the unknown parameter $\rho.$

\begin{proposition}
\label{ident}When the inference on $\rho$ is based on the moment equation
(\ref{mom}), identifiability holds if and only if $\ker\Gamma=\left\{
0\right\}  $.
\end{proposition}

The proof of the Proposition is simple. Let us give a sketch of it now. Assume
that $\ker\Gamma\neq\left\{  0\right\}  $ and pick $v\in\ker\Gamma.$ Setting
$\rho_{v,u}=\rho+v\otimes u$ where $u$ is any vector in $\mathcal{H}$ it is
basic to see that $\Delta=\rho_{v,u}\Gamma$ again. In other words the moment
equation may not be able to distinguish between $\rho$ and $\rho_{v,u}$.

\begin{remark}
The condition $\ker\Gamma=\left\{  0\right\}  $ implies that all the
eigenvalues are strictly positive. In the sequel we will assume that
$\lambda_{1}\geq\lambda_{2}\geq...>0.$
\end{remark}

\subsection{Regularizing the inverse covariance operator}

Even if the identifiability of $\rho$ is ensured by assumption $\mathbf{A}%
_{0}$, we must remain cautious when building an estimator. Several serious
problems appear.

First it is crucial to note that we \textbf{cannot} deduce from (\ref{mom})
that $\Delta\Gamma^{-1}=\rho$. Indeed $\Gamma^{-1}$ does not necessarily
exist. A necessary and sufficient condition for $\Gamma^{-1}$ to be defined as
a linear mapping is : $\ker\Gamma=\left\{  0\right\}  .$ Then $\Gamma^{-1}$ is
an unbounded symmetric operator on $\mathcal{H}$. The consequences are the
following :

\begin{itemize}
\item $\Gamma^{-1}$ is just defined on the dense vector space%
\[
\mathcal{D}\left(  \Gamma^{-1}\right)  =\mathrm{Im}\Gamma=\left\{
x\in\mathcal{H},\ \sum_{i=1}^{n}\dfrac{x_{p}^{2}}{\lambda_{p}^{2}}%
<+\infty\right\}
\]
and $\mathcal{D}\left(  \Gamma^{-1}\right)  \varsubsetneq\mathcal{H}$.

\item $\Gamma^{-1}$ is a measurable linear mapping but is not continuous, in
other words it is continuous at no point for which is it defined or "the
domain of $\Gamma^{-1}$ is also the set of its discontinuities".

\item $\Gamma\Gamma^{-1}$ is not the identity operator on $\mathcal{H}$ but on
$\mathcal{D}\left(  \Gamma^{-1}\right)  $ which entails that (\ref{mom})
implies $\Delta\Gamma^{-1}=\rho_{|\mathrm{Im}\Gamma}\neq\rho$
\end{itemize}

The previous facts are very well-known in operator theory and give rise here
to an ill-posed problem (or an inverse problem). Since $\Gamma^{-1}$ is
extremely irregular, we should propose a way to regularize it i.e. find out
$\Gamma^{\dag}$ say, a linear operator "close" to $\Gamma^{-1}$ and having
additional continuity properties. There are several ways to deal with this
problem. We refer to Arsenin and Tikhonov \cite{ArsTik} and Groetsch
\cite{gro}, amongst many others, for famous books about this topic.

Here the approach is quite intuitive and classical : when (\ref{schmidt})
holds,%
\[
\Gamma^{-1}\left(  x\right)  =\sum_{l\in\mathbb{N}}\dfrac{1}{\lambda_{l}}%
\pi_{l}\left(  x\right)
\]
for all $x$ in $\mathcal{D}\left(  \Gamma^{-1}\right)  .$ We just set
\[
\Gamma^{\dag}\left(  x\right)  =\sum_{l\leq k_{n}}\dfrac{1}{\lambda_{l}}%
\pi_{l}\left(  x\right)
\]
where $\left(  k_{n}\right)  _{n\in N}$ is an increasing sequence tending to
infinity. It may be proved that whenever $x\in\mathcal{D}\left(  \Gamma
^{-1}\right)  $ and $n\uparrow+\infty$,%
\[
\Gamma^{\dag}\left(  x\right)  \rightarrow\Gamma^{-1}\left(  x\right)  .
\]
Besides $\Gamma^{\dag}$ is a continuous operator with $\left\Vert \Gamma
^{\dag}\right\Vert _{\infty}=\lambda_{k_{n}}^{-1}$ and implicitely depends on
$n$.\bigskip

If (\ref{mom}) is the starting point in our estimation procedure, replacing
the unknown operators by their empirical counterparts gives :%
\[
\Delta_{n}=\rho_{n}^{imp}\Gamma_{n}%
\]
where%
\begin{align*}
\Gamma_{n}  &  =\dfrac{1}{n}\sum_{k=1}^{n}X_{k}\otimes X_{k},\\
\Delta_{n}  &  =\dfrac{1}{n-1}\sum_{k=1}^{n-1}X_{k}\otimes X_{k+1}%
\end{align*}
and $\rho_{n}^{imp}$ just \textbf{implicitely} defines our estimate for
$\rho.$

The preceding remarks give some clues to reach the end of the estimation step.
Setting%
\begin{equation}
\Gamma_{n}^{\dag}=\sum_{l\leq k_{n}}\dfrac{1}{\widehat{\lambda}_{l}}%
\widehat{\pi}_{l} \label{gamma_n_dag}%
\end{equation}
where $\widehat{\lambda}_{l}$ and $\widehat{\pi}_{l}$ are the empirical
couterparts of $\lambda_{l}$ and $\pi_{l}$ we get :

\begin{definition}
\label{D}The estimate of $\rho$ is $\rho_{n}$ given by $\rho_{n}=\Delta
_{n}\Gamma_{n}^{\dag}$.
\end{definition}

For further purpose we denote $\widehat{\Pi}_{k_{n}}=\sum_{j=1}^{k_{n}%
}\widehat{\pi}_{j}$ the projector on the space spanned by the $k_{n}$ first
eigenvectors of $\Gamma_{n}$.

\begin{remark}
The $\widehat{\lambda}_{l}$'s and the $\widehat{e}_{l}$'s are obtained as
by-products of the functional PCA of the sample $(X_{1},...,X_{n})$.
\end{remark}

\subsection{A smoothness condition on the autocorrelation operator}

In order to get the main results given in the next section we need to develop
one of the crucial assumptions needed further. This subsection is devoted to
explaining it. This condition must be understood as a smoothness condition on
the unknown operator $\rho$. But what do we mean by "smoothness" for a linear
operator ? The notion of smoothness is intuitively related to functions or
mapping and should be made more clear in our setting. In order to be more
illustrative let us consider for $\rho$ a diagonal operator on $\mathcal{H}$.
Say in any complete othonormal system :%
\[
\rho=diag\left[  \left(  \mu_{i}\right)  _{i\geq1}\right]
\]
with $\mu_{i}\geq\mu_{i+1}.$ Obviously if $\mu_{i}=1$ $\rho=I$ and if the
sequence $\mu_{i}$ is bounded $\left(  \left(  \mu_{i}\right)  _{i\geq1}\in
l^{\infty}\right)  ,$ $\rho$ is a bounded operator. If $\left(  \mu
_{i}\right)  _{i\geq1}\in c_{0},$ $\rho$ is a compact operator. If $\left(
\mu_{i}\right)  _{i\geq1}\in l^{2},$ $\rho$ is a Hilbert-Schmidt operator,
etc. The degree of smoothness of $\rho$ will be strictly determined by the
rate of decrease to zero of $\left(  \left\vert \mu_{i}\right\vert \right)
_{i\geq1}$ or, generally speaking of its eigenvalues or characteristic
numbers. When the $\mu_{i}$'s decrease quickly $\rho$ is "close" to any finite
dimensional approximation based on the $n$ first $\mu_{i}$'s (when $n$ gets
large). Conversely imagine that the $\left\vert \mu_{i}\right\vert $'s tend to
infinity, then $\rho$ is unbounded hence not continuous hence not smooth.

The next assumption
\begin{equation}
\mathbf{A}_{1}:\left\Vert \Gamma^{-1/2}\rho\right\Vert _{\infty}%
<+\infty\label{ass1}%
\end{equation}

tells us that $\rho$ should be at least as "smooth" as $\Gamma^{1/2}.$ Indeed
let us try to be more illustrative and assume that $\rho$ is symmetric and has
the same basis of eigenvectors as $\Gamma$. Assumption (\ref{ass1}) implies
that the sequence $\left(  \mu_{i}/\sqrt{\lambda_{i}}\right)  _{i\in
\mathbb{N}}$ is bounded. We set $\widetilde{\rho}=\Gamma^{-1/2}\rho.$

\begin{remark}
\label{R1}As a consequence of the above we remark for further purpose that if
$\widetilde{\rho}$ is bounded, so is $\widetilde{\rho}^{\ast}.$ But for the
reasons mentioned in the previous subsection $\widetilde{\rho}^{\ast}\neq
\rho^{\ast}\Gamma^{-1/2}$. In fact $\rho^{\ast}\Gamma^{-1/2}$ is a bounded
operator defined on $\mathcal{D}\left(  \Gamma^{-1/2}\right)  $. Like any
bounded operator on a dense domain it may be uniquely extended to a bounded
operator defined on the whole $\mathcal{H}$. This operator precisely coincides
with $\widetilde{\rho}^{\ast}$. I just point out the following : from
(\ref{ass1}) we deduce that%
\begin{equation}
\sup_{p}\left\Vert \rho^{\ast}\Gamma^{-1/2}\left(  e_{p}\right)  \right\Vert
^{2}=\sup_{p}\dfrac{\left\Vert \rho^{\ast}\left(  e_{p}\right)  \right\Vert
^{2}}{\lambda_{p}}\leq M=\left\Vert \widetilde{\rho}^{\ast}\right\Vert
_{\infty}. \label{carrouf}%
\end{equation}

\end{remark}

\section{Main results}

The main results of this work are collected in two theorems below. We first
recapitulate three seminal assumptions under the same label :%
\[
\mathbf{A}_{\mathbf{0}}:\left\{
\begin{array}
[c]{c}%
\ker\Gamma=\left\{  0\right\} \\
\mathbb{E}\left\Vert \varepsilon\right\Vert ^{2}<+\infty\\
\left\Vert \rho\right\Vert _{\infty}<1.
\end{array}
\right.
\]

The subscript $0$ was given on purpose since this set of assumptions is
minimal in order to begin any statistical inference on the model.

Then I remind the reader the so-called Karhunen-Lo\`{e}ve (KL) extension of
the random element $X$ : the distribution of $X$ (i.e. of $X_{n}$ for all $n$
since the sequence is strictly stationary) is :%
\begin{equation}
X=_{d}\sum_{k=1}^{+\infty}\xi_{k}\sqrt{\lambda_{k}}e_{k} \label{KL}%
\end{equation}
where $=_{d}$ denotes equality of distributions and the $\xi_{k}$'s are non
correlated real valued random variables with null expectation and unit
variance (the $\xi_{k}$'s are i.i.d. gaussian if $X_{1}$ is). We will make use
of (\ref{KL}) within the proofs.

The following moment assumption is mild :%
\begin{equation}
\mathbf{A}_{2}:\sup_{k}\mathbb{E}\xi_{k}^{4}<M \label{ass2}%
\end{equation}
It is fullfilled by large families of r.v. $\xi_{k}$'s (subject to
$\mathbb{E}\xi_{k}=0$ and $\mathbb{E}\xi_{k}^{2}=1$) with thin enough queues :
gaussian, uniform, two sided exponential, etc, but will fail for certain
classes of two sided Pareto random variables for instance. Remember that we
study weak convergence for $\rho_{n}\left(  X_{n+1}\right)  ,$ that $\rho_{n}$
depends on $\Gamma_{n}$ and consequently that assumptions on functionals of
the fourth moment of $X_{1}$ (like $\mathbf{A}_{2}$) are unavoidable.

The next assumption is related to the eigenvalues of $\Gamma.$

Let $\lambda_{j}=\lambda\left(  j\right)  $ where $\lambda$ is a positive
function defined on and with values in $\mathbb{R}^{+}$. Clearly function
$\lambda$ is decreasing if the eigenvalues are ordered decreasingly and
$\lim_{t\rightarrow+\infty}\lambda\left(  t\right)  =0$. We assume that :

\begin{center}
$\mathbf{A}_{3}:$\textbf{ The function }$\lambda$\textbf{ is convex}
\end{center}

\begin{remark}
Actually we just need $\mathbf{A}_{3}$ to hold for large values of $j$. This
assumption is finally not constraining at all since it is suited to many
classical cases : when the rate of decay to zero is arithmetic (say
$\lambda_{j}=Const/j^{1+\alpha}$, $\alpha>0$) or exponential ( $\lambda
_{j}=Const\cdot\exp\left(  -\alpha j\right)  $, $\alpha>0$) and in several
other less standard situations such as Laurent series $(\lambda_{j}%
=Const/\left(  j^{\alpha}\log^{1+\beta}j\right)  $, $\alpha,\beta>0$).
\end{remark}

\begin{remark}
Assumption $\mathbf{A}_{3}$ implies that $\lambda_{j}-\lambda_{j+1}\leq
\lambda_{j-1}-\lambda_{j}$.
\end{remark}

The next and first theorem assesses that :

\begin{theorem}
\label{TH1}It is impossible for $\widehat{\rho}_{n}-\rho$ to converge in
distribution for the norm topology on $\mathcal{K}$.
\end{theorem}

\begin{remark}
What is actually proved is : for any normalizing sequence $\alpha_{n}%
\uparrow+\infty,$ $\alpha_{n}\left(  \widehat{\rho}_{n}-\rho\right)  $ either
diverges or converges in distribution to the Dirac distribution on the null
element in $\mathcal{K}$. Also note that weak convergence cannot take place
for the Hilbert-Schmidt topology either since the embedding from
$\mathcal{K}_{2}$ to $\mathcal{K}$ is continuous.
\end{remark}

For technical reasons, we will focus on a sligthly modified version of the
prediction problem. We will assume that $\rho_{n}$ is built from $\left(
X_{1},...,X_{n}\right)  $ and that $X_{n+2}$ is to be predicted from $\rho
_{n}$ and $X_{n+1}.$ In other word the sample is tiled, the last observed
curve ($X_{n+1}$ here) is taken into account to predict $X_{n+2}$ but not to
construct $\rho_{n}$.\newline Here is the main result of the paper. Remind
that $\widehat{\Pi}_{k_{n}}$ was introduced just before Definition \ref{D}.

\begin{theorem}
\label{TH2}When assumptions $\mathbf{A}_{0}-\mathbf{A}_{3}$ hold and if
$k_{n}=o\left(  \dfrac{n^{1/4}}{\log n}\right)  ,$%
\[
\sqrt{\dfrac{n}{k_{n}}}\left(  \widehat{\rho}_{n}\left(  X_{n+1}\right)
-\rho\widehat{\Pi}_{k_{n}}\left(  X_{n+1}\right)  \right)  \overset
{w}{\rightarrow}\mathcal{G}%
\]
where $\mathcal{G}$ is a $\mathcal{H}$-valued gaussian centered random
variable with covariance operator $\Gamma_{\varepsilon}.$
\end{theorem}

\begin{remark}
Theorem \ref{TH1} remains unchanged if $\rho$ is changed into $\rho
\widehat{\Pi}_{k_{n}}$ which appears more "natural" in view of Theorem
\ref{TH2}.
\end{remark}

This central result should be commented. First of all the normalizing sequence
is typically nonparametric : $\sqrt{n/k_{n}}$. Second a bias term appears.
Recently, Cardot, Mas and Sarda \cite{CMS2005} obtained a similar result in a
much simpler regression model, based on i.i.d. observations unlike here. A non
random bias was obtained -namely the random projector $\widehat{\Pi}_{k_{n}}$
was replaced by a non random one- but this could not be carried out here. Also
note that since $\varepsilon$ is the innovation of process $X,$ the best
target we can hope to reach is $\rho\left(  X_{n+1}\right)  $ i.e. the
conditional expectation of $y_{n+1},$ which is random in any case. However it
is simple to prove that $\left\Vert \rho\widehat{\Pi}_{k_{n}}\left(
X_{n+1}\right)  -\rho\left(  X_{n+1}\right)  \right\Vert $ tends to zero in
probability when $n$ tends to infinity. Finally even if the random term
$\rho\widehat{\Pi}_{k_{n}}\left(  X_{n+1}\right)  $ is not quite satisfactory
on a theoretical viewpoint, it may be easily interpreted by practitioners
since $\widehat{\Pi}_{k_{n}}\left(  X_{n+1}\right)  $ is the projection of the
new input onto the $k_{n}$ first axes of the functional PCA of the sample.
These axes have optimality properties w.r.t. the decomposition of variance for
the process $X$.

\section{Concluding remarks}

As seen from the literature on the subject, two modes of stochastic
convergence had already been investigated for estimates of $\rho$ in model
(\ref{modele}) : convergence in probability and almost sure convergence. Weak
convergence was the missing one essentially because it is more
intricate.\newline In fact from%
\[
\left\Vert \rho_{n}\left(  X_{n+1}\right)  -\rho\left(  X_{n+1}\right)
\right\Vert \leq\left\Vert \rho_{n}-\rho\right\Vert _{\infty}\left\Vert
X_{n+1}\right\Vert
\]
it is plain that convergence (almost sure or in probability) for $\left\Vert
\rho_{n}-\rho\right\Vert _{\infty}$ implies convergence for the predictor.
Theorem \ref{TH1} proves that the situation is much more different as far as
convergence in distribution is adressed.\newline It should be also stressed
that assumptions $\mathbf{A}_{0}-\mathbf{A}_{3}$ are truly mild. For instance
all theoretical articles dealing with the problem of asymptotics for the
predictor assume that $\rho$ is symmetric and that the rate of decay of the
sequence of eigenvalues is known.\newline The main advance relies undoubtedly
on the fact that the dimension sequence $k_{n}$ does not depend anymore on the
eigenvalues (previously such conditions as $n^{\alpha}\lambda_{k_{n}%
}\rightarrow+\infty$ for some $\alpha$ where necessary). The existence of a
universal $k_{n}$ enables to revisit all previous results on the topic and
sheds a new light on this model. Indeed in view of Theorem \ref{TH2}, it is
tempting to postulate that a $L^{2}$ minimax rate of convergence could be
$k_{n}/n$ when $\rho$ belongs to the set defined by assumption $\mathbf{A}%
_{1}$ (this set is nothing but an ellipso\"{\i}d of $\mathcal{K}$). But these
considerations are beyond the scope of this article.

\section{Mathematical derivations}

Assumptions $\mathbf{A}_{0}-\mathbf{A}_{3}$ are supposed to hold throughout
the proofs. The generic notation $M$ will be used to denote universal
constants.\ The next equation is straightforward from (\ref{modele}), links
$\Gamma,$ $\Gamma_{\varepsilon}$ and $\rho,$ and will soon be needed :%
\begin{equation}
\Gamma=\rho\Gamma\rho^{\ast}+\Gamma_{\varepsilon}. \label{cov}%
\end{equation}

We start with letting%
\[
S_{n}=\sum_{k=1}^{n}X_{k-1}\otimes\varepsilon_{k}.
\]
Easy calculations give%
\begin{align}
\rho_{n}  &  =\Delta_{n}\Gamma_{n}^{\dagger}=\rho\Gamma_{n}\Gamma_{n}%
^{\dagger}+S_{n}\Gamma_{n}^{\dagger},\nonumber\\
\rho_{n}  &  =\rho\widehat{\Pi}_{k_{n}}+S_{n}\Gamma_{n}^{\dagger}.
\end{align}

It is plain by (\ref{gamma_n_dag}) that $\Gamma_{n}\Gamma_{n}^{\dagger
}=\widehat{\Pi}_{k_{n}}$. Hence :%

\begin{equation}
\rho_{n}-\rho\widehat{\Pi}_{k_{n}}=S_{n}\left(  \Gamma_{n}^{\dagger}%
-\Gamma^{\dagger}\right)  +S_{n}\Gamma^{\dagger} \label{decomp2}%
\end{equation}
which is the starting point.\newline This section is decomposed into three
subsections. In the first one preliminary results and tools connected with the
theory of perturbation for operators on Hilbert spaces are provided. In the
second part I prove that $S_{n}\left(  \Gamma_{n}^{\dagger}-\Gamma^{\dagger
}\right)  $ is a vanishing term if the dimension sequence $k_{n}$ is well
chosen. The third part is devoted to studying weak convergence and proving
Theorem \ref{TH2}. The proof of Theorem \ref{TH1} is postponed to the end of
the paper.

\subsection{Peliminary results}

\subsubsection{Some inequalities}

We first deal with a crucial Lemma.

\begin{lemma}
\label{RA}We have :%
\begin{align}
\sup_{m,p}n\dfrac{\mathbb{E}\left\langle \left(  \Gamma_{n}-\Gamma\right)
\left(  e_{p}\right)  ,e_{m}\right\rangle ^{2}}{\lambda_{p}\lambda_{m}}  &
\leq M\label{RA1}\\
\sup_{m,p}n\dfrac{\mathbb{E}\left\langle S_{n}\left(  e_{p}\right)
,e_{m}\right\rangle ^{2}}{\lambda_{p}\lambda_{m}}  &  \leq M \label{RA2}%
\end{align}

\end{lemma}

\begin{proof}
We begin with proving (\ref{RA1}).%
\begin{align*}
\left\langle \left(  \Gamma_{n}-\Gamma\right)  \left(  e_{p}\right)
,e_{m}\right\rangle ^{2}  &  =\dfrac{1}{n^{2}}\left(  \sum_{k=1}%
^{n}\left\langle X_{k},e_{p}\right\rangle \left\langle X_{k},e_{m}%
\right\rangle \right)  ^{2}\\
\mathbb{E}\left\langle \left(  \Gamma_{n}-\Gamma\right)  \left(  e_{p}\right)
,e_{m}\right\rangle ^{2}  &  =\dfrac{1}{n}\mathbb{E}\left(  \left\langle
X_{1},e_{p}\right\rangle ^{2}\left\langle X_{1},e_{m}\right\rangle ^{2}\right)
\\
&  +\dfrac{2}{n^{2}}\mathbb{E}\sum_{1\leq i<k\leq n}\left(  \left\langle
X_{i},e_{p}\right\rangle \left\langle X_{i},e_{m}\right\rangle \left\langle
X_{k},e_{p}\right\rangle \left\langle X_{k},e_{m}\right\rangle \right)
\end{align*}

It is easily seen by KL decomposition (\ref{KL}) and assumption $\mathbf{A}%
_{2}$ that the first term may be bounded by%
\begin{equation}
\dfrac{1}{n}\lambda_{p}\lambda_{m}\mathbb{E}\left(  \mathbb{\xi}_{p}%
^{2}\mathbb{\xi}_{m}^{2}\right)  \leq M\dfrac{\lambda_{p}\lambda_{m}}{n}
\label{A1}%
\end{equation}
whenever $p=m$ or $p\neq m$.

Now assume that $p\neq m$. We study the second :%
\begin{align*}
X_{k}  &  =\varepsilon_{k}+...+\rho^{k-i-1}\left(  \varepsilon_{i+1}\right)
+\rho^{k-i}\left(  X_{i}\right) \\
X_{k}  &  =E_{k,i}+\rho^{k-i}\left(  X_{i}\right)
\end{align*}
where
\[
E_{k,i}=\varepsilon_{k}+...+\rho^{k-i-1}\left(  \varepsilon_{i+1}\right)
\]
hence%
\begin{align}
&  \mathbb{E}\sum_{i<k}\left(  \left\langle X_{i},e_{p}\right\rangle
\left\langle X_{i},e_{m}\right\rangle \left\langle X_{k},e_{p}\right\rangle
\left\langle X_{k},e_{m}\right\rangle \right) \nonumber\\
&  =\mathbb{E}\sum_{i<k}\left(  \left\langle X_{i},e_{p}\right\rangle
\left\langle X_{i},e_{m}\right\rangle \left\langle \rho^{k-i}\left(
X_{i}\right)  ,e_{p}\right\rangle \left\langle \rho^{k-i}\left(  X_{i}\right)
,e_{m}\right\rangle \right) \nonumber\\
&  +\mathbb{E}\sum_{i<k}\left(  \left\langle X_{i},e_{p}\right\rangle
\left\langle X_{i},e_{m}\right\rangle \left\langle E_{k,i},e_{p}\right\rangle
\left\langle E_{k,i},e_{m}\right\rangle \right) \tag{i}\\
&  =\mathbb{E}\sum_{i<k}\left(  \left\langle X_{i},e_{p}\right\rangle
\left\langle X_{i},e_{m}\right\rangle \left\langle \rho^{k-i}\left(
X_{i}\right)  ,e_{p}\right\rangle \left\langle \rho^{k-i}\left(  X_{i}\right)
,e_{m}\right\rangle \right) \tag{ii}\\
&  =\mathbb{E}\sum_{i<k}\left(  \left\langle X_{1},e_{p}\right\rangle
\left\langle X_{1},e_{m}\right\rangle \left\langle \rho^{k-i}\left(
X_{1}\right)  ,e_{p}\right\rangle \left\langle \rho^{k-i}\left(  X_{1}\right)
,e_{m}\right\rangle \right) \tag{iii}\\
&  =\mathbb{E}\left[  \left\langle X_{1},e_{p}\right\rangle \left\langle
X_{1},e_{m}\right\rangle \sum_{i<k}\left\langle \rho^{k-i}\left(
X_{1}\right)  ,e_{p}\right\rangle \left\langle \rho^{k-i}\left(  X_{1}\right)
,e_{m}\right\rangle \right] \nonumber\\
&  =\mathbb{E}\left[  \left\langle X_{1},e_{p}\right\rangle \left\langle
X_{1},e_{m}\right\rangle \sum_{k=1}^{n-1}\left(  n-k\right)  \left\langle
\rho^{k}\left(  X_{1}\right)  ,e_{p}\right\rangle \left\langle \rho^{k}\left(
X_{1}\right)  ,e_{m}\right\rangle \right]  \tag{iv}%
\end{align}
where (ii) stems from(i) because if $p\neq m$%
\begin{align*}
&  \mathbb{E}\left(  \left\langle X_{i},e_{p}\right\rangle \left\langle
X_{i},e_{m}\right\rangle \left\langle E_{k,i},e_{p}\right\rangle \left\langle
E_{k,i},e_{m}\right\rangle \right) \\
&  =\mathbb{E}\left(  \left\langle X_{i},e_{p}\right\rangle \left\langle
X_{i},e_{m}\right\rangle \right)  \mathbb{E}\left(  \left\langle E_{k,i}%
,e_{p}\right\rangle \left\langle E_{k,i},e_{m}\right\rangle \right) \\
&  =0.
\end{align*}
and (iii) stems from (ii) by stationarity. Now by (iv),%
\begin{align}
&  \dfrac{1}{n}\left\vert \mathbb{E}\sum_{1\leq i<k\leq n}\left(  \left\langle
X_{i},e_{p}\right\rangle \left\langle X_{i},e_{m}\right\rangle \left\langle
X_{k},e_{p}\right\rangle \left\langle X_{k},e_{m}\right\rangle \right)
\right\vert \nonumber\\
&  \leq\mathbb{E}\left[  \left\vert \left\langle X_{1},e_{p}\right\rangle
\left\langle X_{1},e_{m}\right\rangle \right\vert \sum_{k=1}^{n-1}\left(
1-\dfrac{k}{n}\right)  \left\vert \left\langle \rho^{k}\left(  X_{1}\right)
,e_{p}\right\rangle \left\langle \rho^{k}\left(  X_{1}\right)  ,e_{m}%
\right\rangle \right\vert \right]  \label{A2}%
\end{align}
Let us fix $k\geq1$ and develop%
\begin{align*}
\left\langle \rho^{k}\left(  X_{1}\right)  ,e_{p}\right\rangle \left\langle
\rho^{k}\left(  X_{1}\right)  ,e_{m}\right\rangle  &  =\sqrt{\lambda
_{p}\lambda_{m}}\left\langle \Gamma^{-1/2}\rho^{k}\left(  X_{1}\right)
,e_{p}\right\rangle \left\langle \Gamma^{-1/2}\rho^{k}\left(  X_{1}\right)
,e_{m}\right\rangle \\
&  =\sqrt{\lambda_{p}\lambda_{m}}\left\langle \rho^{k-1}\left(  X_{1}\right)
,\dfrac{\rho^{\ast}\left(  e_{p}\right)  }{\sqrt{\lambda_{p}}}\right\rangle
\left\langle \rho^{k-1}\left(  X_{1}\right)  ,\dfrac{\rho^{\ast}\left(
e_{m}\right)  }{\sqrt{\lambda_{m}}}\right\rangle
\end{align*}
and denoting $u_{p}=\rho^{\ast}\left(  e_{p}\right)  /\sqrt{\lambda_{p}}$,%
\begin{align}
\left\vert \left\langle \rho^{k}\left(  X_{1}\right)  ,e_{p}\right\rangle
\left\langle \rho^{k}\left(  X_{1}\right)  ,e_{m}\right\rangle \right\vert  &
\leq\sqrt{\lambda_{p}\lambda_{m}}\left\Vert \rho^{k-1}\right\Vert
^{2}\left\Vert X_{1}\right\Vert ^{2}\left\Vert u_{p}\right\Vert \left\Vert
u_{m}\right\Vert \nonumber\\
&  \leq M\sqrt{\lambda_{p}\lambda_{m}}\left\Vert \rho^{k-1}\right\Vert
^{2}\left\Vert X_{1}\right\Vert ^{2} \label{A3}%
\end{align}
since by (\ref{carrouf}) $\left\Vert u_{p}\right\Vert $ and $\left\Vert
u_{m}\right\Vert $ may be bounded uniformly wrt $p$ and $m$ by $\left\Vert
\mathrm{ext}\left(  \widetilde{\rho}^{\ast}\right)  \right\Vert _{\infty}$
(see Remark \ref{R1}.below) Then%
\begin{align*}
&  \dfrac{1}{n}\left\vert \mathbb{E}\sum_{1\leq i<k\leq n}\left(  \left\langle
X_{i},e_{p}\right\rangle \left\langle X_{i},e_{m}\right\rangle \left\langle
X_{k},e_{p}\right\rangle \left\langle X_{k},e_{m}\right\rangle \right)
\right\vert \\
&  \leq M\sqrt{\lambda_{p}\lambda_{m}}\mathbb{E}\left(  \left\Vert
X_{1}\right\Vert ^{2}\left\vert \left\langle X_{1},e_{p}\right\rangle
\left\langle X_{1},e_{m}\right\rangle \right\vert \right)  \sum_{k=1}%
^{n-1}\left(  1-\dfrac{k}{n}\right)  \left\Vert \rho^{k-1}\right\Vert ^{2}%
\end{align*}
And%
\begin{equation}
\mathbb{E}\left(  \left\Vert X_{1}\right\Vert ^{2}\left\vert \left\langle
X_{1},e_{p}\right\rangle \left\langle X_{1},e_{m}\right\rangle \right\vert
\right)  =\sqrt{\lambda_{p}\lambda_{m}}\left(  \sum_{l=1}^{+\infty}\lambda
_{l}\mathbb{E}\left(  \xi_{l}^{2}\xi_{p}\xi_{m}\right)  \right)  \label{A4}%
\end{equation}
by (\ref{KL}) again. Applying twice Cauchy-Schwarz inequality we bound the
infinite sum by a constant which does not depend on $p$ and $m.$ Collecting
(\ref{A1}), (\ref{A2}), (\ref{A3}) and (\ref{A4}) we get%
\[
n\sup_{p\neq m}\dfrac{\mathbb{E}\left\langle \left(  \Gamma_{n}-\Gamma\right)
\left(  e_{p}\right)  ,e_{m}\right\rangle ^{2}}{\lambda_{p}\lambda_{m}}\leq M
\]
In order to complete the proof (remember that we assumed that $p\neq m$ just
below (\ref{A1})) we can check that our computations remain valid if we take
$p=m$.

The proof of (\ref{RA2}) is similar but simpler. We have%
\begin{align*}
\mathbb{E}\left\langle S_{n}\left(  e_{p}\right)  ,e_{m}\right\rangle ^{2}  &
=\dfrac{1}{n}\mathbb{E}\left(  \left\langle X_{1},e_{p}\right\rangle
^{2}\left\langle \varepsilon_{2},e_{m}\right\rangle ^{2}\right) \\
&  =\dfrac{1}{n}\mathbb{E}\left(  \left\langle X_{1},e_{p}\right\rangle
^{2}\right)  \mathbb{E}\left(  \left\langle \varepsilon_{2},e_{m}\right\rangle
^{2}\right)  =\dfrac{\lambda_{p}\left\langle \Gamma_{\varepsilon}e_{m}%
,e_{m}\right\rangle }{n}\\
&  =\dfrac{\lambda_{p}\left(  \lambda_{m}-\left\langle \rho\Gamma\rho^{\ast
}e_{m},e_{m}\right\rangle \right)  }{n}\leq\dfrac{\lambda_{p}\lambda_{m}}{n}%
\end{align*}

since $\Gamma=\rho\Gamma\rho^{\ast}+\Gamma_{\varepsilon}.$
\end{proof}

The proof of the three following Lemmas may be found in Cardot, Mas, Sarda
\cite{CMS2005}.

\begin{lemma}
\label{trick1}Consider two positive integers $j$ and $k$ large enough and such
that $k>j$. Then
\begin{equation}
j\lambda_{j}\ \geq\ k\lambda_{k}\qquad\mbox{and}\qquad\lambda_{j}-\lambda
_{k}\geq\left(  1-\dfrac{j}{k}\right)  \lambda_{j}. \label{t1}%
\end{equation}
Besides
\begin{equation}
\sum_{j\geq k}\lambda_{j}\leq\left(  k+1\right)  \lambda_{k} . \label{t2}%
\end{equation}

\end{lemma}

\begin{lemma}
\label{trick2}The following is true for $j$ large enough
\[
\sum_{l\neq j}\dfrac{\lambda_{l}}{\left\vert \lambda_{l}-\lambda
_{j}\right\vert }\leq Mj\log j.
\]

\end{lemma}

\subsubsection{A few basic facts about perturbation theory}

Perturbation theory for bounded operators is a powerful tool all along our
study and is of much help when dealing with random (or not) covariance
operators. It features several theoretical interests : for instance
eigenprojectors or pseudo inverses of $\Gamma$ may be expressed as functions
of $\Gamma$ only (without introducing the eigenvectors). However this theory
is not widely used in statistics although the only mathematical prerequisite
is the theory of holomorphic functions and of integrals on contours in the
complex plane. We refer to Dunford-Schwartz \cite{DS} (Chapter VII.3) or to
Gohberg, Goldberg and Kaashoek \cite{GGK} for an introduction to functional
calculus for operators related with Riesz integrals.

Let us denote by $\mathcal{B}_{i}$ the oriented circle of the complex plane
with center $\lambda_{i}$ and radius $\delta_{i}/3$ and define
\[
\mathcal{C}_{n}=\bigcup_{i=1}^{k_{n}}\mathcal{B}_{i}\ .
\]
The open domain whose boundary is $\mathcal{C}_{n}$ is not connected but
however we can apply the functional calculus for bounded operators (see
Dunford-Schwartz \cite{DS}, Section VII.3 Definitions 8 and 9). Results from
perturbation theory yield :%

\[
\Pi_{k_{n}}=\dfrac{1}{2\pi\iota}\int_{\mathcal{C}_{n}}\left(  zI-\Gamma
\right)  ^{-1}dz.
\]
where $\iota^{2}=-1,$ $\Pi_{k_{n}}$ is defined similarly to $\widehat{\Pi
}_{k_{n}}$ (see Theorem \ref{TH2}) and stands for the projector on the space
spanned by the $k_{n}$ first eigenvectors of $\Gamma$. The integral is defined
on the complex plane. Note that the random couterparts (\textit{i.e.}where
$\Pi_{k_{n}}$ and $\Gamma$ are respectively replaced by $\widehat{\Pi}_{k_{n}%
}$ and $\Gamma_{n}$) of the previous equation is just :%
\[
\widehat{\Pi}_{k_{n}}=\dfrac{1}{2\pi\iota}\int_{\widehat{\mathcal{C}}_{n}%
}\left(  zI-\Gamma_{n}\right)  ^{-1}dz.
\]
and the contour $\widehat{\mathcal{C}}_{n}$ is random and depends on the
$\widehat{\lambda}_{i}$'s. The following equalities are also valid%
\begin{align*}
\Gamma^{\dag}  &  =\int_{\mathcal{C}_{n}}z^{-1}\left(  zI-\Gamma\right)
^{-1}dz=\sum_{j=1}^{k_{n}}\int_{\mathcal{B}_{_{j}}}z^{-1}\left(
zI-\Gamma\right)  ^{-1}dz.\\
\Gamma_{n}^{\dag}  &  =\int_{\widehat{\mathcal{C}}_{n}}z^{-1}\left(
zI-\Gamma_{n}\right)  ^{-1}dz=\sum_{j=1}^{k_{n}}\int_{\widehat{\mathcal{B}%
}_{_{j}}}z^{-1}\left(  zI-\Gamma_{n}\right)  ^{-1}dz.
\end{align*}
and%
\begin{align}
&  S_{n}\left(  \Gamma_{n}^{\dagger}-\Gamma^{\dagger}\right)  \left(
X_{n+1}\right) \nonumber\\
&  =\int_{\widehat{\mathcal{C}}_{n}}z^{-1}S_{n}\left(  z-\Gamma_{n}\right)
^{-1}\left(  X_{n+1}\right)  dz-\int_{\mathcal{C}_{n}}z^{-1}S_{n}\left(
z-\Gamma\right)  ^{-1}\left(  X_{n+1}\right)  dz \label{steph}%
\end{align}

As announced at the beginning of the proof section we will prove in the next
subsection that (\ref{steph}) -correctly normalized by $\sqrt{n/k_{n}}$- tends
to zero in probability, hence is negligible. We need two Lemmas to start. In
these Lemmas the square root of a symmetric operator $T,$ say $T^{1/2}$
appears. The bounded operator $T^{1/2}$ has the same eigenvectors as $T.$ Its
eigenvalues are the complex square roots of those of $T.$

\begin{lemma}
\label{thal} We have for $j$ large enough
\begin{align}
\mathbb{E}\sup_{z\in\mathcal{B}_{j}}\left\Vert \left(  zI-\Gamma\right)
^{-1/2}\left(  \Gamma_{n}-\Gamma\right)  \left(  zI-\Gamma\right)
^{-1/2}\right\Vert _{\infty}^{2}  &  \leq\dfrac{M}{n}\left(  j\log j\right)
^{2},\label{majcov}\\
\mathbb{E}\sup_{z\in\mathcal{B}_{j}}\left\Vert z^{-1/2}S_{n}\left(
zI-\Gamma\right)  ^{-1/2}\right\Vert ^{2}  &  \leq\dfrac{M}{n}j\log
j\label{majcroscov}\\
\mathbb{E}\sup_{z\in\mathcal{B}_{j}}\left\Vert \left(  zI-\Gamma\right)
^{-1/2}\varepsilon_{1}\right\Vert ^{2}  &  \leq Mj\log j \label{majeps}%
\end{align}

\end{lemma}

In fact this last Lemma was proved in Cardot, Mas, Sarda \cite{CMS2005} in an
i.i.d framework. However a quick inspection of the proof shows that, by Lemma
\ref{RA} the same result holds in this dependent setting for (\ref{majcov})
and (\ref{majcroscov}). In order to convince the suspicious reader I give now
the derivation of (\ref{majeps}) which uses basically the same technique as
for (\ref{majcov}) and (\ref{majcroscov}) but is shorter. We have :%
\begin{align*}
\left\Vert \left(  zI-\Gamma\right)  ^{-1/2}\varepsilon_{1}\right\Vert ^{2}
&  =\sum_{p=1,p\neq j}^{+\infty}\dfrac{\left\langle \varepsilon_{1}%
,e_{p}\right\rangle ^{2}}{\left\vert z-\lambda_{p}\right\vert }\\
\sup_{z\in\mathcal{B}_{j}}\left\Vert \left(  zI-\Gamma\right)  ^{-1/2}%
\varepsilon_{1}\right\Vert ^{2}  &  =\sum_{p=1,p\neq j}^{+\infty}%
\dfrac{\left\langle \varepsilon_{1},e_{p}\right\rangle ^{2}}{\left\vert
\lambda_{j}-\lambda_{p}\right\vert }%
\end{align*}
since obvioulsy for all $p\neq j,$ $\left\vert z-\lambda_{p}\right\vert
\geq\left\vert \lambda_{j}-\lambda_{p}\right\vert $ when $z\in\mathcal{B}_{j}%
$. Then
\[
\mathbb{E}\sup_{z\in\mathcal{B}_{j}}\left\Vert \left(  zI-\Gamma\right)
^{-1/2}\varepsilon_{1}\right\Vert ^{2}=\sum_{p=1,p\neq j}^{+\infty}%
\dfrac{\mathbb{E}\left\langle \varepsilon_{1},e_{p}\right\rangle ^{2}%
}{\left\vert \lambda_{j}-\lambda_{p}\right\vert }%
\]
Now from $\Gamma=\Gamma_{\varepsilon}+\rho\Gamma\rho^{\ast}$ we see that
$\mathbb{E}\left\langle \varepsilon_{1},e_{p}\right\rangle ^{2}=\left\langle
\Gamma_{\varepsilon}e_{p},e_{p}\right\rangle \leq\left\langle \Gamma
e_{p},e_{p}\right\rangle =\lambda_{p}$ hence%
\[
\mathbb{E}\sup_{z\in\mathcal{B}_{j}}\left\Vert \left(  zI-\Gamma\right)
^{-1/2}\varepsilon_{1}\right\Vert ^{2}\leq\sum_{p=1,p\neq j}^{+\infty}%
\dfrac{\lambda_{p}}{\left\vert \lambda_{j}-\lambda_{p}\right\vert }\leq Mj\log
j
\]
by Lemma \ref{trick2}.

This last Lemma will be used when dealing with residual terms $S_{n}\left(
\Gamma_{n}^{\dagger}-\Gamma^{\dagger}\right)  $ appearing in (\ref{decomp2}).

\begin{lemma}
\label{you}Denoting%
\[
\mathcal{E}_{j}=\left\{  \sup_{z\in\mathcal{B}_{j}}\left\Vert \left(
zI-\Gamma\right)  ^{-1/2}\left(  \Gamma_{n}-\Gamma\right)  \left(
zI-\Gamma\right)  ^{-1/2}\right\Vert _{\infty}<1/2,\right\}  ,
\]
The following holds
\[
\sup_{z\in\mathcal{B}_{j}}\left\Vert \left(  zI-\Gamma\right)  ^{1/2}\left(
zI-\Gamma_{n}\right)  ^{-1}\left(  zI-\Gamma\right)  ^{1/2}\right\Vert
_{\infty}1\!\!1_{\mathcal{E}_{j}}\leq2,\quad a.s.
\]
where $M$ is some positive constant. Besides
\begin{equation}
\mathbb{P}\left(  \mathcal{E}_{j}^{c}\right)  \leq M\dfrac{j\log j}{\sqrt{n}}.
\label{anytime}%
\end{equation}

\end{lemma}

\begin{proof}
We have successively
\[
\left(  zI-\Gamma_{n}\right)  ^{-1}=\left(  zI-\Gamma\right)  ^{-1}+\left(
zI-\Gamma\right)  ^{-1}\left(  \Gamma-\Gamma_{n}\right)  \left(  zI-\Gamma
_{n}\right)  ^{-1},
\]
hence
\[
\left(  zI-\Gamma\right)  ^{1/2}\left(  zI-\Gamma_{n}\right)  ^{-1}\left(
zI-\Gamma\right)  ^{1/2}=I+\left(  zI-\Gamma\right)  ^{-1/2}\left(
\Gamma-\Gamma_{n}\right)  \left(  zI-\Gamma_{n}\right)  ^{-1}\left(
zI-\Gamma\right)  ^{1/2},
\]
and
\begin{equation}
\left[  I+\left(  zI-\Gamma\right)  ^{-1/2}\left(  \Gamma_{n}-\Gamma\right)
\left(  zI-\Gamma\right)  ^{-1/2}\right]  \left(  zI-\Gamma\right)
^{1/2}\left(  zI-\Gamma_{n}\right)  ^{-1}\left(  zI-\Gamma\right)  ^{1/2}=I.
\label{invborn}%
\end{equation}
It is a well known fact that if the linear operator $T$ satisfies $\left\Vert
T\right\Vert _{\infty}<1$ then $I+T$ is an invertible, its inverse is given by
formula
\[
\left(  I+T\right)  ^{-1}=I-T+T^{2}-...
\]
and
\[
\left\Vert \left(  I+T\right)  ^{-1}\right\Vert _{\infty}\leq\dfrac
{1}{1-\left\Vert T\right\Vert _{\infty}}%
\]
From (\ref{invborn}) we deduce that
\begin{align*}
&  \left\Vert \left(  zI-\Gamma\right)  ^{1/2}\left(  zI-\Gamma_{n}\right)
^{-1}\left(  zI-\Gamma\right)  ^{1/2}\right\Vert _{\infty}1\!\!1_{\mathcal{E}%
_{j}}\\
&  =\left\Vert \left[  I+\left(  zI-\Gamma\right)  ^{-1/2}\left(  \Gamma
_{n}-\Gamma\right)  \left(  zI-\Gamma\right)  ^{-1/2}\right]  ^{-1}\right\Vert
_{\infty}1\!\!1_{\mathcal{E}_{j}}\\
&  \leq\dfrac{1}{1-\left\Vert \left(  zI-\Gamma\right)  ^{-1/2}\left(
\Gamma_{n}-\Gamma\right)  \left(  zI-\Gamma\right)  ^{-1/2}\right\Vert
_{\infty}}1\!\!1_{\mathcal{E}_{j}}\leq2,\quad a.s.
\end{align*}
Now, the bound in (\ref{anytime}) stems easily from Markov inequality and
(\ref{majcov}) in Lemma \ref{thal}. This finishes the proof of the Lemma.
\end{proof}

\subsection{Residual term}

This first lemma only aims at proving that the random contour $\widehat
{\mathcal{C}}_{n}$ can be replaced by the non random one $\mathcal{C}_{n}$ in
(\ref{steph}) in order to merge both integrals.

\begin{lemma}
\label{restric}When $\dfrac{1}{\sqrt{n}}k_{n}^{2}\log k_{n}\rightarrow0,$%
\[
S_{n}\left(  \Gamma_{n}^{\dagger}-\Gamma^{\dagger}\right)  \left(
X_{n+1}\right)  =\int_{\mathcal{C}_{n}}z^{-1}S_{n}\left[  \left(  z-\Gamma
_{n}\right)  ^{-1}-\left(  z-\Gamma\right)  ^{-1}\right]  \left(
X_{n+1}\right)  dz+L_{n}%
\]
where $\sqrt{n}\left\Vert L_{n}\right\Vert $ vanishes in probability.
\end{lemma}

\begin{proof}
We introduce the following event :%
\[
\mathcal{A}_{n}=\left\{  \forall j\in\left\{  1,...,k_{n}\right\}
|\dfrac{\left\vert \widehat{\lambda}_{j}-\lambda_{j}\right\vert }{\delta_{j}%
}<1/8\right\}  ,
\]
and $1\!\!1_{\mathcal{A}_{n}}$ is the indicator function of the set
$\mathcal{A}_{n}.$

Introducing the set $\mathcal{A}_{n}$ enables to consider the situation when
all the ordered eigenvalues of $\Gamma_{n}$ are close enough to those of
$\Gamma.$ In fact when $\mathcal{A}_{n}$ holds all the $k_{n}$ first empirical
eigenvalues $\widehat{\lambda}_{j}$ lie in the circle of center $\lambda_{j}$
and radius $\delta_{j}/8,$ say $\widetilde{\mathcal{B}}_{j}$ (included in
$\mathcal{B}_{j}$). Consequently none of the $\widehat{\lambda}_{j}$ is
located in the annulus between $\mathcal{B}_{j}$ and $\widetilde{\mathcal{B}%
}_{j}$ and when $\mathcal{A}_{n}$ holds $\widehat{\mathcal{C}}_{n}$ may be
replaced by $\mathcal{C}_{n}$. It is clear from previous remarks that%
\begin{align*}
S_{n}\left(  \Gamma_{n}^{\dagger}-\Gamma^{\dagger}\right)  \left(
X_{n+1}\right)   &  =S_{n}\left(  \Gamma_{n}^{\dagger}-\Gamma^{\dagger
}\right)  \left(  X_{n+1}\right)  \left(  1\!\!1_{\mathcal{A}_{n}%
}+1\!\!1_{\mathcal{A}_{n}^{c}}\right) \\
&  =\left(  \int_{\mathcal{C}_{n}}z^{-1}S_{n}\left[  \left(  zI-\Gamma
_{n}\right)  ^{-1}-\left(  z-\Gamma\right)  ^{-1}\right]  \left(
X_{n+1}\right)  dz\right) \\
&  -\left(  \int_{\mathcal{C}_{n}}z^{-1}S_{n}\left[  \left(  zI-\Gamma
_{n}\right)  ^{-1}-\left(  z-\Gamma\right)  ^{-1}\right]  \left(
X_{n+1}\right)  dz\right)  1\!\!1_{\mathcal{A}_{n}^{c}}\\
&  +S_{n}\left(  \Gamma_{n}^{\dagger}-\Gamma^{\dagger}\right)  \left(
X_{n+1}\right)  1\!\!1_{\mathcal{A}_{n}^{c}}%
\end{align*}
We set
\begin{align*}
L_{n}  &  =S_{n}\left(  \Gamma_{n}^{\dagger}-\Gamma^{\dagger}\right)  \left(
X_{n+1}\right)  1\!\!1_{\mathcal{A}_{n}^{c}}-\left(  \int_{\mathcal{C}_{n}%
}z^{-1}S_{n}\left[  \left(  zI-\Gamma_{n}\right)  ^{-1}-\left(  zI-\Gamma
\right)  ^{-1}\right]  \left(  X_{n+1}\right)  dz\right)  1\!\!1_{\mathcal{A}%
_{n}^{c}}\\
&  =\left[  S_{n}\Gamma_{n}^{\dagger}\left(  X_{n+1}\right)  -\left(
\int_{\mathcal{C}_{n}}z^{-1}S_{n}\left(  zI-\Gamma_{n}\right)  ^{-1}\left(
X_{n+1}\right)  dz\right)  \right]  1\!\!1_{\mathcal{A}_{n}^{c}}%
\end{align*}
and we see that
\[
\mathbb{P}\left(  \sqrt{n}\left\Vert L_{n}\right\Vert _{\infty}>\varepsilon
\right)  \leq\mathbb{P}\left(  1\!\!1_{\mathcal{A}_{n}^{c}}>\varepsilon
\right)  =\mathbb{P}\left(  \mathcal{A}_{n}^{c}\right)  .
\]
It suffices to get $\mathbb{P}\left(  \mathcal{A}_{n}^{c}\right)
\rightarrow0.$ But
\[
\mathbb{P}\left(  \mathcal{A}_{n}^{c}\right)  \leq\sum_{i=1}^{k_{n}}%
\mathbb{P}\left(  \left\vert \widehat{\lambda}_{i}-\lambda_{i}\right\vert
>\delta_{i}/8\right)  .
\]
Now we refer to Theorem 4.10 of Bosq \cite{Bos3}. Following the proof of this
Theorem along p.122 and 123 it is proved that the asymptotic behaviour of
$\left\vert \widehat{\lambda}_{i}-\lambda_{i}\right\vert $ is the same as
$\left\vert \left\langle \left(  \Gamma_{n}-\Gamma\right)  e_{i}%
,e_{i}\right\rangle \right\vert .$ Then%
\[
\mathbb{P}\left(  \left\vert \widehat{\lambda}_{i}-\lambda_{i}\right\vert
>\delta_{i}/8\right)  \leq8\dfrac{\lambda_{i}}{\delta_{i}}\mathbb{E}\left(
\dfrac{\left\vert \widehat{\lambda}_{i}-\lambda_{i}\right\vert }{\lambda_{i}%
}\right)  \sim8\dfrac{\lambda_{i}}{\delta_{i}}\mathbb{E}\dfrac{\left\vert
\left\langle \left(  \Gamma_{n}-\Gamma\right)  e_{i},e_{i}\right\rangle
\right\vert }{\lambda_{i}}%
\]
By assumption $\mathbf{A}_{2}$ we get%
\[
\mathbb{E}\dfrac{\left\vert \left\langle \left(  \Gamma_{n}-\Gamma\right)
e_{i},e_{i}\right\rangle \right\vert }{\lambda_{i}}\leq\sqrt{\mathbb{E}%
\dfrac{\left\vert \left\langle \left(  \Gamma_{n}-\Gamma\right)  e_{i}%
,e_{i}\right\rangle \right\vert ^{2}}{\lambda_{i}^{2}}}\leq\dfrac{M}{\sqrt{n}}%
\]
by (\ref{RA1}). At last
\[
\mathbb{P}\left(  \mathcal{A}_{n}^{c}\right)  \leq8\dfrac{M}{\sqrt{n}}%
\sum_{i=1}^{k_{n}}\dfrac{\lambda_{i}}{\delta_{i}}\leq\dfrac{M^{\prime}}%
{\sqrt{n}}\sum_{i=1}^{k_{n}}i\log i\leq\dfrac{M^{\prime\prime}}{\sqrt{n}}%
k_{n}^{2}\log k_{n}.
\]

This concludes the proof of the lemma.
\end{proof}

\vspace{.5cm}

\noindent For the sake of clarity, from now on we will abusively note
\[
S_{n}\left(  \Gamma_{n}^{\dagger}-\Gamma^{\dagger}\right)  \left(
X_{n+1}\right)  =\int_{\mathcal{C}_{n}}z^{-1}S_{n}\left[  \left(  z-\Gamma
_{n}\right)  ^{-1}-\left(  z-\Gamma\right)  ^{-1}\right]  \left(
X_{n+1}\right)  dz
\]
but Lemma \ref{restric} above shows that this does not change anything to the
validity of our forthcoming results.

The next Proposition is the central result of this subsection.

\begin{proposition}
\label{bias1}If $\dfrac{1}{\sqrt{n}}k_{n}^{2}\left(  \log k_{n}\right)  ^{2}$
$\rightarrow0$ (which is true if $k_{n}=o\left(  \dfrac{n^{1/4}}{\log
n}\right)  $) we have :%
\[
\sqrt{\dfrac{n}{k_{n}}}S_{n}\left(  \Gamma_{n}^{\dagger}-\Gamma^{\dagger
}\right)  \left(  X_{n+1}\right)  \overset{\mathbb{P}}{\rightarrow}0
\]
in $\mathcal{H}$.
\end{proposition}

\textbf{Proof of Proposition \ref{bias1} :}

We develop :%
\begin{align*}
S_{n}\left(  \Gamma_{n}^{\dagger}-\Gamma^{\dagger}\right)  \left(
X_{n+1}\right)   &  =\int_{\mathcal{C}_{n}}z^{-1}S_{n}\left[  \left(
zI-\Gamma_{n}\right)  ^{-1}-\left(  zI-\Gamma\right)  ^{-1}\right]  \left(
X_{n+1}\right)  dz\\
&  =\int_{\mathcal{C}_{n}}z^{-1}S_{n}\left(  zI-\Gamma\right)  ^{-1}\left(
\Gamma-\Gamma_{n}\right)  \left(  zI-\Gamma_{n}\right)  ^{-1}\left(
X_{n+1}\right)  dz\\
&  =\int_{\mathcal{C}_{n}}z^{-1}S_{n}\left(  zI-\Gamma\right)  ^{-1}\left(
\Gamma-\Gamma_{n}\right)  \left(  zI-\Gamma\right)  ^{-1/2}\\
&  \times\left(  zI-\Gamma\right)  ^{1/2}\left(  zI-\Gamma_{n}\right)
^{-1}\left(  zI-\Gamma\right)  ^{1/2}\left(  zI-\Gamma\right)  ^{-1/2}\left(
X_{n+1}\right)  dz
\end{align*}

and
\begin{align*}
&  \left\Vert S_{n}\left(  \Gamma_{n}^{\dagger}-\Gamma^{\dagger}\right)
\left(  X_{n+1}\right)  \right\Vert \\
&  \leq\int_{\mathcal{C}_{n}}\left\vert z^{-1/2}\right\vert \left\Vert
z^{-1/2}S_{n}\left(  zI-\Gamma\right)  ^{-1/2}\right\Vert _{\infty}\left\Vert
\left(  zI-\Gamma\right)  ^{-1/2}\left(  \Gamma-\Gamma_{n}\right)  \left(
zI-\Gamma\right)  ^{-1/2}\right\Vert _{\infty}\\
&  \times\left\Vert \left(  zI-\Gamma\right)  ^{1/2}\left(  zI-\Gamma
_{n}\right)  ^{-1}\left(  zI-\Gamma\right)  ^{1/2}\right\Vert _{\infty
}\left\Vert \left(  zI-\Gamma\right)  ^{-1/2}\left(  X_{n+1}\right)
\right\Vert dz.
\end{align*}
By Lemma (\ref{you}),%
\begin{align}
&  \left\Vert S_{n}\left(  \Gamma_{n}^{\dagger}-\Gamma^{\dagger}\right)
\left(  X_{n+1}\right)  \right\Vert \nonumber\\
&  =\left\Vert S_{n}\left(  \Gamma_{n}^{\dagger}-\Gamma^{\dagger}\right)
\left(  X_{n+1}\right)  \right\Vert 1\!\!1\left\{  _{\cap_{j}\mathcal{E}_{j}%
}\right\}  +\left\Vert S_{n}\left(  \Gamma_{n}^{\dagger}-\Gamma^{\dagger
}\right)  \left(  X_{n+1}\right)  \right\Vert 1\!\!1_{\left\{  \cup
_{j}\mathcal{E}_{j}^{c}\right\}  }\nonumber\\
&  \leq2\int_{\mathcal{C}_{n}}\left\vert z^{-1/2}\right\vert \left\Vert
z^{-1}S_{n}\left(  zI-\Gamma\right)  ^{-1/2}\right\Vert _{\infty}\left\Vert
\left(  zI-\Gamma\right)  ^{-1/2}\left(  \Gamma-\Gamma_{n}\right)  \left(
zI-\Gamma\right)  ^{-1/2}\right\Vert _{\infty}\label{kafka}\\
&  \times\left\Vert \left(  zI-\Gamma\right)  ^{-1/2}\left(  X_{n+1}\right)
\right\Vert dz+\left\Vert S_{n}\left(  \Gamma_{n}^{\dagger}-\Gamma^{\dagger
}\right)  \left(  X_{n+1}\right)  \right\Vert 1\!\!1_{\cup_{j}\mathcal{E}%
_{j}^{c}}\nonumber
\end{align}
Obviously $\sqrt{n}\left\Vert S_{n}\left(  \Gamma_{n}^{\dagger}-\Gamma
^{\dagger}\right)  \left(  X_{n+1}\right)  \right\Vert 1\!\!1_{\cup
_{j}\mathcal{E}_{j}^{c}}$ decays to zero in probability whenever $\sum
_{j=1}^{k_{n}}\mathbb{P}\left(  \mathcal{E}_{j}^{c}\right)  \rightarrow0$ i.e.
when $\dfrac{1}{\sqrt{n}}\sum_{j=1}^{k_{n}}\left(  j\log j\right)  \sim
\dfrac{k_{n}^{2}\log k_{n}}{\sqrt{n}}$ does.

Let us turn to (\ref{kafka}), tile it into two terms by decomposing $X_{n+1}$
:%
\begin{align*}
W_{1}  &  =\sum_{j=1}^{k_{n}}\int_{\mathcal{B}_{j}}\left\vert z^{-1/2}%
\right\vert \left\Vert z^{-1}S_{n}\left(  zI-\Gamma\right)  ^{-1/2}\right\Vert
_{\infty}\\
&  \times\left\Vert \left(  zI-\Gamma\right)  ^{-1/2}\left(  \Gamma-\Gamma
_{n}\right)  \left(  zI-\Gamma\right)  ^{-1/2}\right\Vert _{\infty}\left\Vert
\left(  zI-\Gamma\right)  ^{-1/2}\left(  \varepsilon_{n+1}\right)  \right\Vert
dz\\
W_{2}  &  =\sum_{j=1}^{k_{n}}\int_{\mathcal{B}_{j}}\left\vert z^{-1/2}%
\right\vert \left\Vert z^{-1}S_{n}\left(  zI-\Gamma\right)  ^{-1/2}\right\Vert
_{\infty}\\
&  \times\left\Vert \left(  zI-\Gamma\right)  ^{-1/2}\left(  \Gamma-\Gamma
_{n}\right)  \left(  zI-\Gamma\right)  ^{-1/2}\right\Vert _{\infty}\left\Vert
\left(  zI-\Gamma\right)  ^{-1/2}\rho\left(  X_{n}\right)  \right\Vert dz
\end{align*}
and first prove that $\sqrt{n/k_{n}}W_{1}$ tends in probability to zero. Let
us simplifiy this first term.%
\begin{align*}
W_{1}  &  \leq\sum_{j=1}^{k_{n}}\dfrac{\delta_{j}}{\sqrt{\left\vert
\lambda_{j}-\delta_{j}\right\vert }}\sup_{z\in\mathcal{B}_{j}}\left\{
\left\Vert z^{-1/2}S_{n}\left(  zI-\Gamma\right)  ^{-1/2}\right\Vert _{\infty
}\left\Vert \left(  zI-\Gamma\right)  ^{-1/2}\left(  \varepsilon_{n+1}\right)
\right\Vert \right\} \\
&  \times\sup_{z\in\mathcal{B}_{j}}\left\{  \left\Vert \left(  zI-\Gamma
\right)  ^{-1/2}\left(  \Gamma-\Gamma_{n}\right)  \left(  zI-\Gamma\right)
^{-1/2}\right\Vert _{\infty}\right\}
\end{align*}
hence%
\begin{align}
\mathbb{E}W_{1}  &  \leq\sum_{j=1}^{k_{n}}\dfrac{\delta_{j}}{\sqrt{\left\vert
\lambda_{j}-\delta_{j}\right\vert }}\sqrt{\mathbb{E}\sup_{z\in\mathcal{B}_{j}%
}\left\{  \left\Vert \left(  zI-\Gamma\right)  ^{-1/2}\left(  \Gamma
-\Gamma_{n}\right)  \left(  zI-\Gamma\right)  ^{-1/2}\right\Vert _{\infty
}\right\}  ^{2}}\nonumber\\
&  \times\sqrt{\mathbb{E}\sup_{z\in\mathcal{B}_{j}}\left\{  \left\Vert
z^{-1/2}S_{n}\left(  zI-\Gamma\right)  ^{-1/2}\right\Vert _{\infty}\right\}
^{2}\mathbb{E}\sup_{z\in\mathcal{B}_{j}}\left\{  \left\Vert \left(
zI-\Gamma\right)  ^{-1/2}\left(  \varepsilon_{n+1}\right)  \right\Vert
\right\}  ^{2}}\tag{i}\\
&  \leq\dfrac{M}{n}\sum_{j=1}^{k_{n}}\dfrac{\delta_{j}}{\sqrt{\left\vert
\lambda_{j}-\delta_{j}\right\vert }}\left(  j\log j\right)  ^{2}\tag{ii}\\
&  \leq\dfrac{M}{n}\sum_{j=1}^{k_{n}}\sqrt{\delta_{j}}j^{2}\left(  \log
j\right)  ^{2}\leq\dfrac{M}{n}k_{n}^{5/2}\left(  \log k_{n}\right)
^{2}\nonumber
\end{align}
From (i) to (ii) I invoke Lemma \ref{thal}, $\dfrac{\delta_{j}}{\sqrt
{\left\vert \lambda_{j}-\delta_{j}\right\vert }}$ was bounded by $\sqrt
{\delta_{j}}$, at last it is plain that $\sqrt{j\delta_{j}}$ is bounded. As a
consequence of the above if one chooses $k_{n}$ such that
\[
\sqrt{\dfrac{n}{k_{n}}}\dfrac{1}{n}k_{n}^{5/2}\left(  \log k_{n}\right)
^{2}=\dfrac{1}{\sqrt{n}}k_{n}^{2}\left(  \log k_{n}\right)  ^{2}\rightarrow0
\]
we see that $\sqrt{\dfrac{n}{k_{n}}}W_{1}$ tends in probability to zero. We
turn to the second term $W_{2}$ and like above%
\begin{align*}
W_{2}  &  \leq\sum_{j=1}^{k_{n}}\dfrac{\delta_{j}}{\sqrt{\left\vert
\lambda_{j}-\delta_{j}\right\vert }}\sup_{z\in\mathcal{B}_{j}}\left\{
\left\Vert z^{-1/2}S_{n}\left(  zI-\Gamma\right)  ^{-1/2}\right\Vert _{\infty
}\left\Vert \left(  zI-\Gamma\right)  ^{-1/2}\Gamma^{1/2}\widetilde{\rho
}\left(  X_{n}\right)  \right\Vert \right\} \\
&  \times\sup_{z\in\mathcal{B}_{j}}\left\{  \left\Vert \left(  zI-\Gamma
\right)  ^{-1/2}\left(  \Gamma-\Gamma_{n}\right)  \left(  zI-\Gamma\right)
^{-1/2}\right\Vert _{\infty}\right\}
\end{align*}

The situation is slightly more complicated than above since $S_{n}$ is not
independent from $X_{n}.$ We introduce a truncation. Assume that $\tau_{n}$ is
an increasing sequence tending to infinity.%
\begin{align*}
W_{2}  &  =W_{2}\mathrm{I}_{\left\{  \left\Vert X_{n}\right\Vert <\tau
_{n}\right\}  }+W_{2}\mathrm{I}_{\left\{  \left\Vert X_{n}\right\Vert \geq
\tau_{n}\right\}  }\\
&  =W_{2}^{-}+W_{2}^{+}.
\end{align*}
Obviously $\sqrt{\dfrac{n}{k_{n}}}W_{2}^{+}$ tends in probability to zero
since for all $\varepsilon>0$%
\[
\mathbb{P}\left(  \sqrt{\dfrac{n}{k_{n}}}W_{2}^{+}>\varepsilon\right)
\leq\mathbb{P}\left(  \left\Vert X_{n}\right\Vert \geq\tau_{n}\right)
\leq\dfrac{\mathbb{E}\left\Vert X_{1}\right\Vert }{\tau_{n}}.
\]
We turn to
\begin{align*}
W_{2}^{-}  &  \leq\left\Vert \widetilde{\rho}\right\Vert \tau_{n}\sum
_{j=1}^{k_{n}}\dfrac{\delta_{j}}{\sqrt{\left\vert \lambda_{j}-\delta
_{j}\right\vert }}\sup_{z\in\mathcal{B}_{j}}\left\{  \left\Vert z^{-1/2}%
S_{n}\left(  zI-\Gamma\right)  ^{-1/2}\right\Vert _{\infty}\left\Vert \left(
zI-\Gamma\right)  ^{-1/2}\Gamma^{1/2}\right\Vert _{\infty}\right\} \\
&  \times\sup_{z\in\mathcal{B}_{j}}\left\{  \left\Vert \left(  zI-\Gamma
\right)  ^{-1/2}\left(  \Gamma-\Gamma_{n}\right)  \left(  zI-\Gamma\right)
^{-1/2}\right\Vert _{\infty}\right\}  ,\\
\mathbb{E}W_{2}^{-}  &  \leq\left\Vert \widetilde{\rho}\right\Vert \tau
_{n}\sum_{j=1}^{k_{n}}\dfrac{\delta_{j}}{\sqrt{\left\vert \lambda_{j}%
-\delta_{j}\right\vert }}\sup_{z\in\mathcal{B}_{j}}\left\{  \left\Vert \left(
zI-\Gamma\right)  ^{-1/2}\Gamma^{1/2}\right\Vert _{\infty}\right\}
\sqrt{\mathbb{E}\sup_{z\in\mathcal{B}_{j}}\left\{  \left\Vert z^{-1/2}%
S_{n}\left(  zI-\Gamma\right)  ^{-1/2}\right\Vert _{\infty}^{2}\right\}  }\\
&  \times\sqrt{\mathbb{E}\sup_{z\in\mathcal{B}_{j}}\left\{  \left\Vert \left(
zI-\Gamma\right)  ^{-1/2}\left(  \Gamma-\Gamma_{n}\right)  \left(
zI-\Gamma\right)  ^{-1/2}\right\Vert _{\infty}^{2}\right\}  }\\
&  \leq M\dfrac{\tau_{n}}{n}\sum_{j=1}^{k_{n}}\dfrac{\delta_{j}}%
{\sqrt{\left\vert \lambda_{j}-\delta_{j}\right\vert }}j^{2}\left(  \log
j\right)  ^{3/2}\leq M\dfrac{\tau_{n}k_{n}^{5/2}\left(  \log k_{n}\right)
^{3/2}}{n}%
\end{align*}
hence $\sqrt{\dfrac{n}{k_{n}}}W_{2}^{-}$ tends in probability to zero whenever
$\dfrac{\tau_{n}k_{n}^{2}\left(  \log k_{n}\right)  ^{3/2}}{\sqrt{n}%
}\rightarrow0.$ Now we choose $\tau_{n}=\sqrt{\log k_{n}}$ with $k_{n}$ as
above for $W_{1}$. This finishes the proof of Proposition \ref{bias1}.

\subsection{Weakly convergent term}

As seen from (\ref{decomp2}) and from previous subsection $S_{n}%
\Gamma^{\dagger}\left(  X_{n+1}\right)  $ will fully determine the asymptotics
of the predictor :%
\begin{align*}
S_{n}\Gamma^{\dagger}\left(  X_{n+1}\right)   &  =\sum_{k=1}^{n}\left\langle
X_{k-1},\Gamma^{\dagger}\left(  X_{n+1}\right)  \right\rangle \varepsilon
_{k}\\
&  =\sum_{k=1}^{n}Z_{k,n}.
\end{align*}

We decompose $Z_{k,n}$ in three terms%
\begin{align*}
Z_{k,n}  &  =Z_{k,n}^{+}+Z_{k,n}^{0}+Z_{k,n}^{-}\\
Z_{k,n}^{+}  &  =\left\langle \Gamma^{\dagger}X_{k-1},\varepsilon_{n+1}%
+\rho\left(  \varepsilon_{n}\right)  +...+\rho^{n-k}\left(  \varepsilon
_{k+1}\right)  \right\rangle \varepsilon_{k}\\
Z_{k,n}^{0}  &  =\left\langle \Gamma^{\dagger}X_{k-1},\rho^{n+1-k}%
\varepsilon_{k}\right\rangle \varepsilon_{k}\\
Z_{k,n}^{-}  &  =\left\langle \Gamma^{\dagger}X_{k-1},\rho^{n+2-k}\left(
X_{k-1}\right)  \right\rangle \varepsilon_{k}%
\end{align*}
stemming from
\[
X_{n+1}=\varepsilon_{n+1}+\rho\left(  \varepsilon_{n}\right)  +...+\rho
^{n+1-k}\left(  \varepsilon_{k}\right)  +\rho^{n+2-k}\left(  X_{k-1}\right)
.
\]
We will show in Lemma \ref{remain} below that the series involving
$Z_{k,n}^{0}$ and $Z_{k,n}^{-}$ are negligible ; weak convergence is strictly
determined by $\sum_{k=1}^{n}Z_{k,n}^{+}.$ The asymptotic distribution is
given at Proposition \ref{weakconv} below. We begin with an important Lemma.

\begin{lemma}
\label{martindiff}The random sequences $Z_{k,n}^{+}$ and $Z_{k,n}^{-}$ are
Hilbert-valued martingale difference arrays w.r.t. the sequence $\left(
\mathcal{F}_{i}\right)  _{i\leq k}$ where $\mathcal{F}_{i}$ is the $\sigma
$-algebra generated by $\left(  \varepsilon_{l}\right)  _{l\leq i}$
\end{lemma}

\textbf{Proof :}

Denoting%
\begin{align*}
X_{k,n}^{\sharp}  &  =\varepsilon_{n+1}+\rho\left(  \varepsilon_{n}\right)
+...+\rho^{n-k}\left(  \varepsilon_{k+1}\right)  ,\\
\mathbb{E}\left(  Z_{k,n}^{+}|\mathcal{F}_{k-1}\right)   &  =\mathbb{E}\left(
\left\langle \Gamma^{\dagger}X_{k-1},X_{k,n}^{\sharp}\right\rangle
\varepsilon_{k}|\mathcal{F}_{k-1}\right)  .
\end{align*}
Since $\varepsilon_{k}$ is independent from $X_{k,n}^{\sharp}$ and both
sequences of random elements are centered we deduce that%
\[
\mathbb{E}\left(  Z_{k,n}^{+}|\mathcal{F}_{k-1}\right)  =0.
\]
Then
\begin{align*}
\mathbb{E}\left(  Z_{k,n}^{-}|\mathcal{F}_{k-1}\right)   &  =\mathbb{E}\left(
\left\langle \Gamma^{\dagger}X_{k-1},\rho^{n+2-k}\left(  X_{k-1}\right)
\right\rangle \varepsilon_{k}|\mathcal{F}_{k-1}\right) \\
&  =\left\langle \Gamma^{\dagger}X_{k-1},\rho^{n+2-k}\left(  X_{k-1}\right)
\right\rangle \mathbb{E}\left(  \varepsilon_{k}|\mathcal{F}_{k-1}\right) \\
&  =0
\end{align*}

\begin{proposition}
\label{weakconv}%
\[
S_{n}^{+}=\dfrac{1}{\sqrt{nk_{n}}}\sum_{k=1}^{n}Z_{k,n}^{+}\overset
{w}{\rightarrow}\mathcal{G}\left(  0,\Gamma_{\varepsilon}\right)  .
\]

\end{proposition}

\textbf{Proof of the Proposition :}

Since $\sum_{k=1}^{n}Z_{k,n}^{+}$ is a $\mathcal{H}$-valued martingale
difference array we first could hope to apply existing criteria for weak
convergence of such sequences. Most of these criteria (see Walk \cite{Walk} or
Rackauskas \cite{Ra1}) rely on convergence in probability for the conditional
covariance operator. They do not seem to be adapted in this context (I could
not go through with it...). I propose the reader to come back to the "sources"
of the Central Limit Theorem on infinite dimensional vector spaces. We will
simply prove that $S_{n}^{+}$ is a uniformly tight sequence and that finite
distributions, when computed on a sufficiently large set of functionals
converge to gaussian limits, hence characterizing the limiting covariance
operator $\Gamma_{\varepsilon}$. In order to understand this approach I refer
to the paper by A. de Acosta \cite{Ac}, especially to Theorem 2.3 p.279.

For further purpose we begin with a first Lemma in which covariance and
cross-covariance operators for the array $Z_{k,n}^{+}$ are computed.

\begin{lemma}
\label{covcov}If $k<i,$ $\mathbb{E}\left(  Z_{k,n}^{+}\otimes Z_{i,n}%
^{+}\right)  =0$ and
\[
\mathbb{E}\left(  Z_{k,n}^{+}\otimes Z_{k,n}^{+}\right)  =\Gamma_{\varepsilon
}\left(  k_{n}-\mathrm{tr}\left(  \Gamma^{\dagger}\rho^{n-k+1}\Gamma\left(
\rho^{\ast}\right)  ^{n-k+1}\right)  \right)  .
\]

\end{lemma}

\begin{proof}%
\[
Z_{k,n}^{+}\otimes Z_{i,n}^{+}=\left\langle \Gamma^{\dagger}X_{k-1}%
,X_{k,n}^{\sharp}\right\rangle \left\langle \Gamma^{\dagger}X_{i-1}%
,X_{i,n}^{\sharp}\right\rangle \left(  \varepsilon_{k}\otimes\varepsilon
_{i}\right)
\]
and since $X_{i-1}=\rho^{i-k}\left(  X_{k-1}\right)  +\varepsilon
_{i-1}+...+\rho^{i-1-k}\left(  \varepsilon_{k}\right)  .$ We tile $Z_{k,n}%
^{+}\otimes Z_{i,n}^{+}$ into two terms. We see that%
\[
\mathbb{E}\left[  \left\langle \Gamma^{\dagger}X_{k-1},X_{k,n}^{\sharp
}\right\rangle \left\langle \Gamma^{\dagger}\rho^{i-k}\left(  X_{k-1}\right)
,X_{i,n}^{\sharp}\right\rangle \left(  \varepsilon_{k}\otimes\varepsilon
_{i}\right)  \right]  =0
\]
since $\varepsilon_{k}$ is independent from all the other terms. The second
term is :%
\[
\left\langle \Gamma^{\dagger}X_{k-1},X_{k,n}^{\sharp}\right\rangle
\left\langle \Gamma^{\dagger}\left(  \varepsilon_{i-1}+...+\rho^{i-1-k}\left(
\varepsilon_{k}\right)  \right)  ,X_{i,n}^{\sharp}\right\rangle \left(
\varepsilon_{k}\otimes\varepsilon_{i}\right)  .
\]
Its expectation is null since $X_{k-1}$ is centered and independent from all
the other terms. We focus on the second part of the Lemma.

We have
\begin{align*}
\mathbb{E}\left(  Z_{k,n}^{+}\otimes Z_{k,n}^{+}\right)   &  =\left(
\mathbb{E}\left\langle \Gamma^{\dagger}X_{k-1},X_{k,n}^{\sharp}\right\rangle
^{2}\right)  \mathbb{E}\left(  \varepsilon_{k}\otimes\varepsilon_{k}\right) \\
&  =\left(  \mathbb{E}\left\langle \Gamma^{\dagger}X_{k-1},X_{k,n}^{\sharp
}\right\rangle ^{2}\right)  \Gamma_{\varepsilon}%
\end{align*}

and
\begin{align*}
\mathbb{E}\left\langle \Gamma^{\dagger}X_{k-1},X_{k,n}^{\sharp}\right\rangle
^{2}  &  =\mathbb{E}\left(  \mathbb{E}\left\langle X_{k-1},\Gamma^{\dagger
}X_{k,n}^{\sharp}\right\rangle ^{2}|X_{k,n}^{\sharp}\right) \\
&  =\mathbb{E}\left\Vert \Gamma^{1/2}\Gamma^{\dagger}X_{k,n}^{\sharp
}\right\Vert ^{2}\\
&  =\mathbb{E}\left\Vert \Gamma^{\dagger1/2}X_{k,n}^{\sharp}\right\Vert ^{2}\\
&  =\mathrm{tr}\left(  \Gamma^{\dagger}\Gamma_{k,n}^{\sharp}\right)
\end{align*}
where
\begin{align*}
\Gamma_{k,n}^{\sharp}  &  =\mathbb{E}\left(  X_{k,n}^{\sharp}\otimes
X_{k,n}^{\sharp}\right) \\
&  =\Gamma_{\varepsilon}+\rho\Gamma_{\varepsilon}\rho^{\ast}+...+\rho
^{n-k}\Gamma_{\varepsilon}\left(  \rho^{\ast}\right)  ^{n-k}\\
&  =\Gamma-\rho^{n-k+1}\Gamma\left(  \rho^{\ast}\right)  ^{n-k+1}.
\end{align*}
Then
\begin{align*}
\mathrm{tr}\left(  \Gamma^{\dagger}\Gamma_{k,n}^{\sharp}\right)   &
=\mathrm{tr}\left(  \Gamma^{\dagger}\Gamma\right)  -\mathrm{tr}\left(
\Gamma^{\dagger}\rho^{n-k+1}\Gamma\left(  \rho^{\ast}\right)  ^{n-k+1}\right)
\\
&  =k_{n}-\mathrm{tr}\left(  \Gamma^{\dagger}\rho^{n-k+1}\Gamma\left(
\rho^{\ast}\right)  ^{n-k+1}\right)  .
\end{align*}
The proof of Lemma \ref{covcov} is complete.
\end{proof}

Now we prove that all the finite-dimensional distributions converge to a
gaussian limit. It suffices to get, for all $x$ in $\mathcal{H}$,%

\begin{equation}
\dfrac{1}{\sqrt{nk_{n}}}\sum_{k=1}^{n}\left\langle Z_{k,n}^{+},x\right\rangle
\overset{w}{\rightarrow}\mathcal{N}\left(  0,\sigma_{\varepsilon,x}%
^{2}\right)  \label{wcfindim}%
\end{equation}
where $\sigma_{\varepsilon,x}^{2}=\mathbb{E}\left\langle \varepsilon
_{k},x\right\rangle ^{2}.$

Since $\left\langle Z_{k,n}^{+},x\right\rangle $ is a real valued MDA it
suffices to apply the criteria given in Mac Leish \cite{McL}. In view of Lemma
(\ref{covcov}) it is enough to prove that $\sum_{k=1}^{n}\mathrm{tr}\left(
\Gamma^{\dagger}\Gamma_{k,n}^{\sharp}\right)  \sim nk_{n}$ that is
\[
\dfrac{\sum_{k=1}^{n}\mathrm{tr}\left(  \Gamma^{\dagger}\Gamma_{k,n}^{\sharp
}\right)  -nk_{n}}{nk_{n}}=\dfrac{\sum_{k=1}^{n}\mathrm{tr}\left(
\Gamma^{\dagger}\rho^{n-k+1}\Gamma\left(  \rho^{\ast}\right)  ^{n-k+1}\right)
}{nk_{n}}\rightarrow0
\]
The usual properties of the trace provide%
\begin{align*}
\left\vert \mathrm{tr}\left(  \Gamma^{\dagger}\rho^{n-k+1}\Gamma\left(
\rho^{\ast}\right)  ^{n-k+1}\right)  \right\vert  &  =\left\vert
\mathrm{tr}\left(  \left(  \rho^{\ast}\right)  ^{n-k+1}\Gamma^{\dagger}%
\rho^{n-k+1}\Gamma\right)  \right\vert \\
&  =\left\vert \mathrm{tr}\left(  \left(  \rho^{\ast}\right)  ^{n-k+1}%
\Gamma^{\dagger}\rho^{n-k+1}\Gamma\right)  \right\vert \\
&  \leq\left\Vert \left(  \rho^{\ast}\right)  ^{n-k+1}\Gamma^{\dagger}%
\rho^{n-k+1}\right\Vert _{\infty}\left\vert \mathrm{tr}\Gamma\right\vert \\
&  =\left\Vert \left(  \rho^{\ast}\right)  ^{n-k}\widetilde{\rho}^{\ast}%
\Gamma^{1/2}\Gamma^{\dagger}\Gamma^{1/2}\widetilde{\rho}\rho^{n-k}\right\Vert
_{\infty}\left\vert \mathrm{tr}\Gamma\right\vert \\
&  \leq\left\Vert \rho^{n-k}\right\Vert ^{2}\left\Vert \widetilde{\rho}^{\ast
}\right\Vert _{\infty}\left\Vert \widetilde{\rho}\right\Vert _{\infty
}\left\vert \mathrm{tr}\Gamma\right\vert
\end{align*}
and we see that whenever $nk_{n}\rightarrow+\infty$%
\begin{equation}
\sum_{k=1}^{n}\left(  \mathbb{E}\left\langle \Gamma^{\dagger}X_{k-1}%
,X_{k,n}^{\sharp}\right\rangle ^{2}\right)  \sim nk_{n} \label{nicole}%
\end{equation}
which ensures (\ref{wcfindim}).

Now we turn to the second part of the proof, namely : "the sequence $\left(
S_{n}^{+}\right)  _{n\in\mathbb{N}}$ is tight". Once more we go through a Lemma.

\begin{lemma}
\label{flatconc}By $\mathcal{P}_{m}$ we denote the projector associated to the
$m$ first eigenvectors of the covariance operator $\Gamma_{\varepsilon}$ of
$\varepsilon_{1}.$ Then,%
\begin{equation}
\limsup_{m\rightarrow+\infty}\sup_{n}\mathbb{P}\left(  \left\Vert \left(
I-\mathcal{P}_{m}\right)  S_{n}^{+}\right\Vert >\varepsilon\right)  =0.
\label{convtight}%
\end{equation}

\end{lemma}

\begin{remark}
What we prove is "with prescribed probability the sequence $S_{n}^{+}$ is
concentrated in the $\varepsilon$-neighborhood of a finite dimensional space
-i.e. \textrm{Im}$\left(  \mathcal{P}_{m}\right)  $". This phenomenon is
called flat concentration and ensures the tightness of $\left(  S_{n}%
^{+}\right)  _{n\in\mathbb{N}}$ (see de Acosta (1970), Definition 2.1 p.279).
\end{remark}

\textbf{Proof of Lemma \ref{flatconc}} :%

\[
\mathbb{P}\left(  \left\Vert \left(  I-\mathcal{P}_{m}\right)  S_{n}%
^{+}\right\Vert >\varepsilon\right)  \leq\dfrac{\mathbb{E}\left(  \left\Vert
\left(  I-\mathcal{P}_{m}\right)  S_{n}^{+}\right\Vert ^{2}\right)
}{\varepsilon^{2}}%
\]
where
\begin{align}
\mathbb{E}\left(  \left\Vert \left(  I-\mathcal{P}_{m}\right)  S_{n}%
^{+}\right\Vert ^{2}\right)   &  =\dfrac{1}{nk_{n}}\mathbb{E}\left(
\left\Vert \sum_{k=1}^{n}\left\langle \Gamma^{\dagger}X_{k-1},\varepsilon
_{n+1}+\rho\left(  \varepsilon_{n}\right)  +...+\rho^{n-k}\left(
\varepsilon_{k+1}\right)  \right\rangle \left(  I-\mathcal{P}_{m}\right)
\varepsilon_{k}\right\Vert ^{2}\right) \tag{i}\\
&  =\dfrac{1}{nk_{n}}\left(  \sum_{k=1}^{n}\mathbb{E}\left(  \left\langle
\Gamma^{\dagger}X_{k-1},X_{k,n}^{\sharp}\right\rangle ^{2}\left\Vert \left(
I-\mathcal{P}_{m}\right)  \varepsilon_{k}\right\Vert ^{2}\right)  \right)
\tag{ii}\\
&  =\dfrac{1}{nk_{n}}\mathbb{E}\left\Vert \left(  I-\mathcal{P}_{m}\right)
\varepsilon_{k}\right\Vert ^{2}\left(  \sum_{k=1}^{n}\mathbb{E}\left\langle
\Gamma^{\dagger}X_{k-1},X_{k,n}^{\sharp}\right\rangle ^{2}\right) \nonumber\\
&  =\dfrac{1}{nk_{n}}\mathrm{tr}\left(  \left(  I-\mathcal{P}_{m}\right)
\Gamma_{\varepsilon}\right)  \left(  \sum_{k=1}^{n}\mathbb{E}\left\langle
\Gamma^{\dagger}X_{k-1},X_{k,n}^{\sharp}\right\rangle ^{2}\right)  .\nonumber
\end{align}
On line (ii) the expectation of all the cross products is null. I skip through
these calculations since they are exactly alike thoses carried within Lemma
\ref{covcov} above. The computations made in the first part of the proof (see
display (\ref{nicole})) are useful here. They ensure that
\[
\sup_{n}\dfrac{1}{nk_{n}}\left(  \sum_{k=1}^{n}\mathbb{E}\left\langle
\Gamma^{\dagger}X_{k-1},X_{k,n}^{\sharp}\right\rangle ^{2}\right)  <M
\]
where $M$ is some universal constant. At last letting $m$ tend to infinity we
get%
\[
\lim_{m\rightarrow+\infty}\mathrm{tr}\left(  \left(  I-\mathcal{P}_{m}\right)
\Gamma_{\varepsilon}\right)  =0
\]
which proves Lemma \ref{flatconc}.

It remains to conclude. Lemma \ref{flatconc} ensures that the centered
sequence $S_{n}^{+}$ is tight.\ By (\ref{wcfindim}) we know that the weak
limit is gaussian and that its covariance function (hence its covariance
operator) is fully characterized : the same as $\varepsilon_{1}$. We invoke
for instance A. de Acosta (1970) to conclude the proof of Proposition
\ref{weakconv}.

\begin{lemma}
\label{remain}%
\begin{align}
&  \dfrac{1}{\sqrt{nk_{n}}}\sum_{k=1}^{n}Z_{k,n}^{-}\overset{\mathbb{P}%
}{\rightarrow}0,\label{rest1}\\
&  \dfrac{1}{\sqrt{nk_{n}}}\sum_{k=1}^{n}Z_{k,n}^{0}\overset{\mathbb{P}%
}{\rightarrow}0. \label{rest2}%
\end{align}

\end{lemma}

\textbf{Proof }$:$

It is plain that $Z_{k,n}^{-}$ is an array of non-correlated random elements.
We prove that%
\[
\dfrac{1}{nk_{n}}\mathbb{E}\left\Vert \sum_{k=1}^{n}Z_{k,n}^{-}\right\Vert
^{2}\rightarrow0
\]%
\begin{align*}
\mathbb{E}\left\Vert \sum_{k=1}^{n}Z_{k,n}^{-}\right\Vert ^{2}  &
=\mathbb{E}\left\Vert \varepsilon_{1}\right\Vert ^{2}\sum_{k=1}^{n}%
\mathbb{E}\left\langle \Gamma^{\dagger}X_{k-1},\rho^{n+2-k}\left(
X_{k-1}\right)  \right\rangle ^{2}\\
&  =\mathbb{E}\left\Vert \varepsilon_{1}\right\Vert ^{2}\sum_{k=1}%
^{n}\mathbb{E}\left\langle \left(  \Gamma^{\dagger}\right)  ^{1/2}%
X_{k-1},\Gamma^{-1/2}\rho^{n+2-k}\left(  X_{k-1}\right)  \right\rangle ^{2}\\
&  \leq\mathbb{E}\left\Vert \varepsilon_{1}\right\Vert ^{2}\sum_{k=1}%
^{n}\mathbb{E}\left[  \left\Vert \left(  \Gamma^{\dagger}\right)
^{1/2}X_{k-1}\right\Vert ^{2}\left\Vert X_{k-1}\right\Vert ^{2}\right]
\left\Vert \Gamma^{-1/2}\rho^{n+2-k}\right\Vert _{\infty}\\
&  =\left\Vert \widetilde{\rho}\right\Vert _{\infty}\mathbb{E}\left\Vert
\varepsilon_{1}\right\Vert ^{2}\mathbb{E}\left[  \left\Vert \left(
\Gamma^{\dagger}\right)  ^{1/2}X_{1}\right\Vert ^{2}\left\Vert X_{1}%
\right\Vert ^{2}\right]  \sum_{k=1}^{n}\left\Vert \rho^{n+1-k}\right\Vert
_{\infty}.
\end{align*}
Since KL expansion yields%
\[
\left\Vert \left(  \Gamma^{\dagger}\right)  ^{1/2}X_{1}\right\Vert
^{2}\left\Vert X_{1}\right\Vert ^{2}=_{d}\sum_{i=1}^{k_{n}}\xi_{i}^{2}%
\sum_{j=1}^{+\infty}\lambda_{j}\xi_{j}^{2}%
\]
we easily see by assumption $\mathbf{A}_{2}$ that
\begin{equation}
\mathbb{E}\left[  \left\Vert \left(  \Gamma^{\dagger}\right)  ^{1/2}%
X_{1}\right\Vert ^{2}\left\Vert X_{1}\right\Vert ^{2}\right]  =O\left(
k_{n}\right)  \label{choz}%
\end{equation}
hence (\ref{rest1}).

We turn to obtaining a bound for the second term. With $Z_{k,n}^{0}%
=\left\langle \Gamma^{\dagger}X_{k-1},\rho^{n+1-k}\varepsilon_{k}\right\rangle
\varepsilon_{k}$ we get :%
\begin{align*}
\mathbb{E}\left\Vert \sum_{k=1}^{n}Z_{k,n}^{0}\right\Vert ^{2}  &  =\sum
_{k=1}^{n}\mathbb{E}\left\Vert Z_{k,n}^{0}\right\Vert ^{2}+2\sum_{1\leq
i<j\leq n}\mathbb{E}\left\langle Z_{i,n}^{0},Z_{j,n}^{0}\right\rangle \\
&  =\sum_{k=1}^{n}\mathbb{E}\left\langle \Gamma^{\dagger}X_{k-1},\rho
^{n+1-k}\varepsilon_{k}\right\rangle ^{2}\left\Vert \varepsilon_{k}\right\Vert
^{2}\\
&  +2\sum_{1\leq i<j\leq n}\mathbb{E}\left(  \left\langle \Gamma^{\dagger
}X_{i-1},\rho^{n+1-i}\varepsilon_{i}\right\rangle \left\langle \varepsilon
_{i},\varepsilon_{j}\right\rangle \left\langle \Gamma^{\dagger}X_{j-1}%
,\rho^{n+1-j}\varepsilon_{j}\right\rangle \right)  .
\end{align*}
The first term may be bounded by%
\begin{align*}
&  \sum_{k=1}^{n}\mathbb{E}\left[  \left\Vert \left(  \Gamma^{\dagger}\right)
^{1/2}\rho^{n+1-k}\varepsilon_{k}\right\Vert ^{2}\left\Vert \varepsilon
_{k}\right\Vert ^{2}\right] \\
&  \leq\left\Vert \widetilde{\rho}\right\Vert _{\infty}\mathbb{E}\left(
\left\Vert \varepsilon_{1}\right\Vert ^{4}\right)  \sum_{k=1}^{n}\left\Vert
\rho^{n-k}\right\Vert _{\infty}.
\end{align*}
The second term may be rewritten :%
\begin{align*}
&  \sum_{1\leq i<j\leq n}\mathbb{E}\left\langle \Gamma^{\dagger}X_{i-1}%
,\rho^{n+1-i}\varepsilon_{i}\right\rangle \left\langle \varepsilon
_{i},\varepsilon_{j}\right\rangle \left\langle \Gamma^{\dagger}X_{j-1}%
,\rho^{n+1-j}\varepsilon_{j}\right\rangle \\
&  =\sum_{1\leq i<j\leq n}\mathbb{E}\left\langle \Gamma^{\dagger}X_{i-1}%
,\rho^{n+1-i}\varepsilon_{i}\right\rangle \left\langle \Gamma^{\dagger}%
\rho^{n+1-j}\Gamma_{\varepsilon}\left(  \varepsilon_{i}\right)  ,X_{j-1}%
\right\rangle \\
&  =\sum_{1\leq i<j\leq n}\mathbb{E}\left\langle \Gamma^{\dagger}X_{i-1}%
,\rho^{n+1-i}\varepsilon_{i}\right\rangle \left\langle \Gamma^{\dagger}%
\rho^{n+1-j}\Gamma_{\varepsilon}\left(  \varepsilon_{i}\right)  ,\rho
^{j-i}X_{i-1}\right\rangle \\
&  =\sum_{1\leq i<j\leq n}\mathbb{E}\left\langle \left(  \Gamma^{\dagger
}\right)  ^{1/2}X_{i-1},\widetilde{\rho}\rho^{n-i}\varepsilon_{i}\right\rangle
\left\langle \widetilde{\rho}\rho^{n-j}\Gamma_{\varepsilon}\left(
\varepsilon_{i}\right)  ,\widetilde{\rho}\rho^{j-i-1}X_{i-1}\right\rangle .
\end{align*}
Taking absolute values we get the bound%
\[
\left\Vert \widetilde{\rho}\right\Vert _{\infty}^{3}\mathbb{E}\left(
\left\Vert \varepsilon_{1}\right\Vert ^{2}\right)  \mathbb{E}\left[
\left\Vert X_{1}\right\Vert \left\Vert \left(  \Gamma^{\dagger}\right)
^{1/2}X_{1}\right\Vert \right]  \left\Vert \Gamma_{\varepsilon}\right\Vert
_{\infty}\sum_{1\leq i<j\leq n}\left\Vert \rho^{n-i}\right\Vert _{\infty
}\left\Vert \rho^{n-j}\right\Vert _{\infty}\left\Vert \rho^{j-i-1}\right\Vert
_{\infty}.
\]
Once again invoking (\ref{choz}) we get (\ref{rest2}) in Lemma \ref{remain}%
.\medskip

\textbf{Proof of Theorem \ref{TH1} :}

From all that was done above it is straightforward to deduce that weak
convergence for $\rho_{n}-\rho$ depends only on the term $S_{n}\Gamma
^{\dagger}$ in (\ref{decomp2}). We recall it : $S_{n}\Gamma^{\dagger}%
=\sum_{k=1}^{n}\Gamma^{\dagger}X_{k-1}\otimes\varepsilon_{k}.$ I guess the
reader will agree with the following sentences : "Assume that $\left(
\varepsilon_{k}\right)  _{k\in\mathbb{Z}}$ and $\left(  X_{k}\right)
_{k\in\mathbb{Z}}$ are independent sequences of independent random elements in
$\mathcal{H}$. Then if in this framework $S_{n}\Gamma^{\dagger}$ does not
converge weakly $S_{n}\Gamma^{\dagger}$ will not converge weakly in the
setting of model (\ref{modele})". Obviously the situation is much favourable
assuming independence "everywhere".

Let us assume that $\dfrac{\alpha_{n}}{n}S_{n}\Gamma^{\dagger}$ converges
weakly to some random variable $Z$ for some increasing sequence $\alpha_{n}.$
We deduce that, for any $f\in\mathcal{K}^{\ast},$ the dual space of
$\mathcal{K}^{\ast},$
\[
f\left(  \dfrac{\alpha_{n}}{n}S_{n}\Gamma^{\dagger}\right)  =\dfrac{\alpha
_{n}}{n}f\left(  S_{n}\Gamma^{\dagger}\right)  =\dfrac{\alpha_{n}}{n}%
\sum_{k=1}^{n}f\left(  \Gamma^{\dagger}X_{k-1}\otimes\varepsilon_{k}\right)
\]
converges weakly to $f\left(  Z\right)  .$ In fact $\mathcal{K}^{\ast
}=\mathcal{K}_{1}$ the space of trace class operators (see Dunford-Schwartz
\cite{DS} for this classical result), the duality bracket is nothing than the
usual trace. Consquently we should investigate weak convergence for
\[
\dfrac{\alpha_{n}}{n}\sum_{k=1}^{n}\mathrm{tr}\left[  T\left(  \Gamma
^{\dagger}X_{k-1}\otimes\varepsilon_{k}\right)  \right]  =\dfrac{\alpha_{n}%
}{n}\sum_{k=1}^{n}\left\langle \Gamma^{\dagger}X_{k-1},T\varepsilon
_{k}\right\rangle
\]
where $T$ is a trace class operator. To prove Theorem \ref{TH1}, it is enough
to take $T=u\otimes v,$ $u,v\in\mathcal{H}.$ Indeed%
\[
f\left(  \dfrac{\alpha_{n}}{n}S_{n}\Gamma^{\dagger}\right)  =\dfrac{\alpha
_{n}}{n}\sum_{k=1}^{n}\left\langle \Gamma^{\dagger}X_{k-1},T\varepsilon
_{k}\right\rangle =\dfrac{\alpha_{n}}{n}\sum_{k=1}^{n}\left\langle
X_{k-1},\Gamma^{\dagger}v\right\rangle \left\langle u,\varepsilon
_{k}\right\rangle
\]
Now we consider two cases depending on the location of $v$ :

\begin{enumerate}
\item If $v\in D\left(  \Gamma^{-1}\right)  ,$ $\Gamma^{\dagger}v$ is a
bounded sequence that converges to $\Gamma^{-1}v.$ It is straightforward to
see that $f\left(  \dfrac{1}{\sqrt{n}}S_{n}\Gamma^{\dagger}\right)  $
converges in distribution to $f\left(  Z\right)  $ (which is gaussian) by the
real CLT for i.i.d. r.v. \textbf{This means that necessarily }$\alpha
_{n}=\sqrt{n}$.

\item Let us take a general $v\notin D\left(  \Gamma^{-1}\right)  ,$ and
compute the variance of the series above with $\alpha_{n}=\sqrt{n}$%
\begin{align*}
\mathbb{E}\left[  f\left(  \dfrac{1}{\sqrt{n}}S_{n}\Gamma^{\dagger}\right)
\right]  ^{2} &  =\dfrac{1}{n}\sum_{k=1}^{n}\mathbb{E}\left\langle
X_{k-1},\Gamma^{\dagger}v\right\rangle ^{2}\mathbb{E}\left\langle
u,\varepsilon_{k}\right\rangle ^{2}\\
&  =\dfrac{\sigma_{\varepsilon,u}^{2}}{n}\sum_{k=1}^{n}\mathbb{E}\left\langle
X_{k-1},\Gamma^{\dagger}v\right\rangle ^{2}\\
&  =\sigma_{\varepsilon,u}^{2}\left\Vert \Gamma^{1/2}\Gamma^{\dagger
}v\right\Vert ^{2}=\sigma_{\varepsilon,u}^{2}\left\Vert \left(  \Gamma
^{\dagger}\right)  ^{1/2}v\right\Vert ^{2}%
\end{align*}
where $\sigma_{\varepsilon,u}^{2}=\mathbb{E}\left\langle u,\varepsilon
_{1}\right\rangle ^{2}$ and%
\[
\left\Vert \left(  \Gamma^{\dagger}\right)  ^{1/2}v\right\Vert ^{2}=\sum
_{i=1}^{k_{n}}\dfrac{\left\langle v,e_{i}\right\rangle ^{2}}{\lambda_{i}}%
\]
Choosing $\left\langle v,e_{i}\right\rangle ^{2}=\lambda_{i}$ or $\left\langle
v,e_{i}\right\rangle ^{2}=\lambda_{i}\beta_{i}$ where $\beta_{i}%
\rightarrow\beta>0$ we see that $\left\Vert \left(  \Gamma^{\dagger}\right)
^{1/2}v\right\Vert ^{2}\rightarrow+\infty$ and the real valued random variable
$f\left(  \left(  1/\sqrt{n}\right)  S_{n}\Gamma^{\dagger}\right)  $ cannot
converge weakly since its variance tends to infinity. This shows that the
marginals of $\left(  \alpha_{n}/n\right)  S_{n}\Gamma^{\dagger}$ do not all
converge to the same limiting measure and not all at the same rate, which
prevents weak convergence in the topology of $\mathcal{K}$. Hence Theorem
\ref{TH1}.\bigskip
\end{enumerate}

\end{document}